\documentclass[reqno]{amsart}
\usepackage{bbm}
\usepackage{amsmath,xcolor,esint,geometry}
\usepackage[foot]{amsaddr}
\usepackage{amssymb}
\usepackage{mathtools}
\usepackage{bookmark}
\def\Xint#1{\mathchoice
{\XXint\displaystyle\textstyle{#1}}
{\XXint\textstyle\scriptstyle{#1}}
{\XXint\scriptstyle\scriptscriptstyle{#1}}
{\XXint\scriptscriptstyle\scriptscriptstyle{#1}}
\!\int}
\def\XXint#1#2#3{{\setbox0=\hbox{$#1{#2#3}{\int}$ }
\vcenter{\hbox{$#2#3$ }}\kern-.6\wd0}}

\def\dashint{\Xint-}
\DeclarePairedDelimiter\abs{\lvert}{\rvert}

\mathtoolsset{showonlyrefs}
 \allowdisplaybreaks
\def\R{\mathbb R}
\def\H{\mathcal H}
\def\ecc{\mathbf{e}_C}

\def\Om{\Omega}
\def\N{\mathbb N}
\def\e{\varepsilon}

\newcommand{\D}{\nabla}
\newcommand{\dd}{\partial}
\newcommand{\norm}[1]{\left\lVert#1\right\rVert}

\newtheorem{theorem}{Theorem}
\newtheorem{lemma}{Lemma}
\newtheorem*{theoremm*}{Theorem}
\newtheorem{remark}{Remark}
\newtheorem{definition}{Definition}

\begin{document}

\title{Regularity Results for a free interface problem with H\"{o}lder coefficients}

\author[1]{L. Esposito}
\author[2]{L. Lamberti}

\begin{abstract}
{We study a class of variational problems involving both bulk and interface energies. The bulk energy is of Dirichlet type albeit of very general form allowing the dependence from the unknown variable $u$ and the position $x$. We employ the regularity theory of $\Lambda$-minimizers to study the regularity of the free interface. The hallmark of the paper is the mild regularity assumption concerning the dependence of the coefficients with respect to $x$ and $u$ that is of H\"{o}lder type. }
\end{abstract}

\maketitle

	
\makeatletter

\makeatother

\section{Introduction and statements}
\label{Introduction and statements}
This paper deals with a large class of nonlinear variational problems involving both
bulk and interface energies,
\begin{equation}\label{intro0}
{\mathcal F}(E,u;\Omega)=\int_\Om\bigl [ F(x,u,\nabla u)+\mathbbm{1}_{E}G(x,u,\nabla u)\bigr ]\,\,dx +P(E;\Om)\,,
\end{equation}
where $u\in H^1(\Omega)$ and $\mathbbm{1}_{E}$ denotes the characteristic function of a set $E\subset \Omega$ with finite perimeter 
$P(E;\Omega)$ in $\Omega$. Energy functionals including both bulk and interface terms are very frequent in mathematical and physical literature (see for instance \cite{AC}, \cite{AB}, \cite{FF}, \cite{FFLM}, \cite{Gur}, \cite{Lar}, \cite{LH}, \cite{Lin}, \cite{Tay}). In particular, the functionals that we study in this paper are strictly related to the integral energy employed in the study of charged droplets (see \cite{DHV}, \cite{MV}). A prototype version of these functionals, 
that is
\begin{equation}\label{model}
\int_{\Omega}\sigma_{E}(x)|\nabla u|^2\,dx+ P(E;\Omega),
\end{equation}
with $u=u_0$ prescribed on $\partial \Omega$ and $\sigma_{E}(x)=\beta \mathbbm{1}_{E}+\alpha \mathbbm{1}_{\Omega\setminus E}$, $0<\alpha<\beta$,
 was formerly studied in {1993} in two papers by L. Ambrosio \& G. Buttazzo and F.H. Lin (see \cite{AB} and  \cite{Lin}).\\
The regularity of minimizers of these kinds of functionals is a rather subtle
issue even in the scalar setting especially regarding the free interface $\partial E$.\\
In {1993} in the paper \cite{AB} L. Ambrosio and G. Buttazzo proved that if $(E,u)$ is a minimizer of the functional $\eqref{model}$, then $u$ is locally H\"{o}lder continuous in $\Omega$ and $E$ is relatively open in $\Omega$. In the same volume of the same journal, F.H. Lin proved a regularity result for the interface $\partial E$.
To clarify the situation we define the set of regular points of $\partial E$ as follows:
\begin{equation}
\mbox{Reg}(E):= \left\{x\in \partial E\cap \Omega\,:\,\partial E \text{ is a }C^{1,\gamma} \mbox{ hypersurface in some } I(x) \text{ and for some }\gamma\in(0,1)
\right\}
\end{equation}
where $I(x)$ denotes a neighborhood of $x$. Accordingly, we define the set of singular points of $\partial E$
\begin{equation}
\Sigma(E) := (\partial E \cap \Omega)\setminus \mbox{Reg}(E).
\end{equation}
In \cite{Lin} F.H. Lin proved that, for minimal configurations of the functional $\eqref{model}$, 
$$\mathcal {H}^{n-1}(\Sigma(E))=0.
$$
The aforementioned regularity result has been recently improved by G. De Philippis \& A. Figalli, and N. Fusco \& V. Julin. Using different approaches and different techniques G. De Philippis \& A. Figalli in \cite{DF} and N. Fusco \& V. Julin in \cite{FJ} proved that for minimal configurations of the functional \eqref{model} it turns out that,
\begin{equation}
\label{redu}
dim_{\mathcal{H}}(\Sigma(E))\leq n-1-\varepsilon,
\end{equation}
for some $\varepsilon>0$ depending only on $\alpha, \beta$. Regarding this dependence, it is worth noticing that in \cite{EF} it was proven that $u\in C^{0,\frac{1}{2}+\varepsilon}$ and the reduced boundary $\dd^*E$ of $E$ is a $C^{1,\varepsilon}-$hypersurface and $\mathcal{H}^s(\partial E\setminus \partial^* E)=0$ for all $s>n-8$, assuming that $1\leq\frac{\alpha}{\beta}<\gamma_n$, for some ${\gamma_n}>1$ depending only on the dimension.\\
In 1999 F.H. Lin and R.V. Kohn in \cite{LK} extended the same result that the first author obtained for the model case $\eqref{model}$ to the more general setting of integral energy of the type $\eqref{intro0}$, depending also on $x$ and $u$. More precisely F.H. Lin and R.V. Kohn proved, for minimal configurations $(E,u)$ of \eqref{intro0} under suitable smothness assumption of $F$ and $G$, that $\mathcal {H}^{n-1}(\Sigma(E))=0$. \\
A natural question to ask is whether the same dimension reduction of the singular set $\Sigma(E)$ proved for the model case $\eqref{model}$ by G. De Philippis \& A. Figalli and N. Fusco \& V. Julin can be extended also to the general case of functionals of the type $\eqref{intro0}$.  In a very recent paper we give a positive answer to this question.
Inded in \cite{EL} we prove that
$$
dim_{\mathcal{H}}(\Sigma(E))\leq n-1-\varepsilon,
$$
for some $\varepsilon>0$, for optimal configurations of a wide class of quadratic functionals depending also on $x$ and $u$.
Our path to prove the aforementioned result basically follows the same strategy used in \cite{FJ}. The technique used in \cite{EL} relies on the linearity of the Euler-Lagrange equation of the functional $\eqref{intro0}$. For this reason we need a quadratic structure condition for the bulk energy. Conversely, the nonquadratic case is less studied and there are few regularity results available (see \cite{CFP}, \cite{CFP2}, \cite{E}, \cite{Lam}).

Throughout the paper we will assume that the density energies $F$ and $G$ in $\eqref{intro0}$ satisfy the following structural quadratic assumptions:

\begin{align}\label{structure1}
&F(x,s,z)=\sum_{i,j=1}^n a_{ij}(x,s)z_iz_j+\sum_{i=1}^n a_i(x,s)z_i+a(x,s),\\
\label{structure2}
&G(x,s,z)=\sum_{i,j=1}^n b_{ij}(x,s)z_iz_j+\sum_{i=1}^n b_i(x,s)z_i+b(x,s),
\end{align}
for any $(x,s,z)\in\Omega\times\R\times\R^n$. In the paper \cite{EL} we assumed as in \cite{LK} that the coefficients $a_{ij},b_{ij},,a_i,b_i,a,b$ belong to the class ${C^{0,1}(\Omega\times \R)}$ with respect to both variables $x$ and $s$. This ${C^{0,1}}$ assumption of the coefficients with respect to $(x,s)$ is crucial in several respect in order to prove the desired regularity result for $\partial E$.\\
In the first place the ${C^{0,1}}$ assumption is strongly used (see Theorem 2 in \cite{EL}) to prove that every minimizer of the constrained problem (that is for $|E|=d$ fixed) is a $\Lambda$-minimizer of a penalized functional containing the extraterm $\Lambda||E|-d|$.
In addition the ${C^{0,1}}$ assumption is primarly used to get an Euler-Lagrange-type equations that is one of the main ingredients to prove the desired regularity result (see Proposition 4.9 in \cite{FJ}) and Theorem 8 in \cite{EL}).\\
In this paper we examine in depth the question of the minimal regularity assumptions of the coeficients we ought to assume in order to get the regularity result quoted in \eqref{redu}.  Concerning the coefficients appearing in \eqref{structure1} and \eqref{structure2} we will assume H\"{o}lder continuous dependence of $(x,s)$. We exploited the proof strategy in every possible way in order to push to the limit the assumptions concerning the H\"{o}lder exponent of the coefficients. In this regard it is important to point out that no restriction is needed for the H\"{o}lder exponent $\beta$ with respect to the $s$ variable quoted below.
Precisely we will assume that
\begin{equation*}
a_{ij}(x,\cdot),b_{ij}(x,\cdot),a_i(x,\cdot),b_i(x,\cdot),a(x,\cdot),b(x,\cdot)\in C^{0,\beta}(\R),\;\; \text{ for every }x\in\Omega.
\end{equation*}
We will denote by $L_\beta$ the greatest H\"older seminorm of the coefficients with respect to the second variable, that is
\begin{equation}
\label{Hoelderianity1}
[a_{ij}(x,\cdot)]_{\beta}:=\sup_{u,t\in \R, \,u\neq t}\frac{|a_{ij}(x,u)-a_{ij}(x,t)|}{|u-t|}\leq L_\beta,\quad \forall x\in\Omega,
\end{equation}
and the same holds true for $b_{ij},a_i,b_i,a,b$.\\
Similarly we will assume about the dependence on the first variable,
\begin{equation*}
a_{ij}(\cdot,s),b_{ij}(\cdot,s),a_i(\cdot,s),b_i(\cdot,s),a(\cdot,s),b(\cdot,s)\in C^{0,{\alpha}}(\Omega),\;\; \text{ for every }s\in\R,
\end{equation*} 
where 
$$\alpha\in \left(\frac{n-1}{n},1\right].$$
We will denote by $L_\alpha$ the greatest H\"older seminorm of the coefficients with respect to the first variable, that is
\begin{equation}\label{Hoelderianity2}
[a_{ij}(\cdot,s)]_{\alpha}:=\sup_{y,z\in \Omega,\, y\neq z}\frac{|a_{ij}(y,s)-a_{ij}(z,s)|}{|y-z|}\leq L_\alpha,\quad \forall s\in \R,
\end{equation}
and the same holds true for $b_{ij},a_i,b_i,a,b$.\\
Moreover, to ensure the existence of minimizers we assume the boundedness of the coefficients and the ellipticity of the matrices $a_{ij}$ and $b_{ij}$,
\begin{align}\label{ellipticity1}
&\nu|z|^2\leq a_{ij}(x,s)z_i z_j\leq N |z|^2,\quad\nu |z|^2\leq b_{ij}(x,s)z_i z_j\leq  N |z|^2,\\
\label{ellipticity2}
& \sum_{i=1}^n |a_i(x,s)|+\sum_{i=1}^n |b_i(x,s)|+|a(x,s)|+|b(x,s)|\leq L,
\end{align}
for any $(x,s,z)\in\Omega\times\R\times\R^n$, where $\nu$, $N$ and $L$ are three positive constants. \\
Some comments about the H\"{o}lder exponent $\alpha$ are in order. There are two main points in our proof where the assumption $\alpha>\frac{n-1}{n}$ is used. In both cases we have to handle with a perturbation of the set $E$.\\ 
The first point concerns the equivalence between the constrained problem and the penalized problem  (see the definitions below). In Theorem \ref{Teorema Penalizzazione} we perform a suitable ``small'' perturbation of a minimal set $E$ around a point $x\in \partial E$ using a transformation of the type
$$
\Phi_{\sigma}(x)=x+\sigma X(x), \;\;\;\;\text{ where }X\in C^1_0(B_r(x)).
$$
If we denote by ${\widetilde E}:=\Phi_{\sigma}(E)$ the perturbed set and by $\tilde{u}:=u \circ \Phi_{\sigma}^{-1}$ the perturbed function, we prove that
$$
{\mathcal F}(E,u)-{\mathcal F}({\widetilde E},{\tilde u})=O(\sigma^{\alpha}),
$$
where $\alpha$ is the H\"{o}lder exponent given in \eqref{Hoelderianity2}. On the other hand, in Theorem \ref{Teorema Penalizzazione} we prove by contradiction that $(E,u)$ is a minimizer of a penalized functional obtained adding in \eqref{intro0} a penalization term of the type
$$
\Lambda\big||\tilde{E}| -d\big|^s,
$$
for some suitable $\Lambda$ to be choosen sufficiently large. Since we can observe that  $\Lambda\big||\tilde{E}| -d\big|^s=O(\sigma^{s})$, it is clear that we are forced to choose $s=\sigma$ (see Definition \ref{Penalizzato} below). Finally it is evident that this new penalization term cannot exceed the perimeter term when we rescale the functional (see Lemma \ref{Lemma riscalamento}) and so we are forced to choose $\alpha>\frac{n-1}{n}$.
\\
The second point concerns the excess improvement given in Theorem \ref{Miglioramento eccesso}, where we use a standard rescaling argument to show that the limit $g$ of the rescaled functions whose graph locally represents $\partial{E}$ is armonic (see Step 1 in Theorem \ref{Miglioramento eccesso}). In this step we use the Taylor expansion of the bulk term given in Theorem \ref{VarB} and the condition $\alpha>\frac{n-1}{n}$ is again crucial, see \eqref{5}.\\

In this paper we study the regularity of minimizers of the following constrained problem.
\begin{definition}
\label{ProblemaVincolato} We shall denote by $\eqref{P_c}$ the constrained problem
\begin{equation}
\label{P_c}
\min_{\substack{E\in \mathcal{A}(\Omega)\\v\in u_0+H_0^1(\Omega)}}
\left\{
\mathcal{F}(E,v;\Omega)\,:\,  |E|=d
\right\},
\tag{$P_c$}
\end{equation}
where $u_0\in H^1(\Omega)$, $0<d<|\Omega|$ are given and $\mathcal{A}(\Omega)$ is the class of all subsets of $\Omega$ with finite perimeter in $\Omega$.
\end{definition}
The problem of handling with the constraint $|E|=d$ is overtaken using an argument introduced in \cite{EF}, ensuring that every minimizer of the constrained problem $\eqref{P_c}$ is also a minimizer of a penalized functional of the type
\begin{equation}\label{intro1}
{\mathcal F}_{\Lambda}(E,v;\Omega)=\mathcal{F}(E,v;\Omega)+\Lambda\big||E| -d\big|^\alpha,
\end{equation}
for some suitable $\Lambda>0$ (see Theorem \ref{Teorema Penalizzazione} below). Therefore, we give in addition the following definition.
\begin{definition}\label{Penalizzato}
We shall denote by $\eqref{P}$ the penalized problem
\begin{equation}
\label{P}
\min_{\substack{E\in \mathcal{A}(\Omega)\\v\in u_0+H_0^1(\Omega)}}
\mathcal{F}_{\Lambda}(E,v;\Omega),
\tag{$P$}
\end{equation}
where $u_0\in H^1(\Omega)$ is fixed and $\mathcal{A}(\Omega)$ is the same class defined in Definition \ref{ProblemaVincolato}.
\end{definition}
From the point of view of regularity, the extra term $\Lambda\big||E| -d\big|^{\alpha}$ is a higher order negligible perturbation, being $\alpha>\frac{n-1}{n}$.
The main result of the paper is stated in the following theorem.
\begin{theorem}
\label{Teorema principale}
Let $(E,u)$ be a minimizer of problem $\eqref{P}$, under assumptions $\eqref{structure1}-\eqref{ellipticity2}$. Then
\begin{itemize}
\item[a)] there exists a relatively open set $\Gamma\subset \partial E$ such that $\Gamma$ is a $C^{1,\mu}$ hypersurface for all $0<\mu<\frac{\gamma}{2}$, where 
{$\gamma:=1+n(\alpha-1)\in(0,1)$} ,
\item[b)] there exists $\varepsilon >0$ depending on $n,\nu,N,L$, such that $$\mathcal{H}^{n-1-\varepsilon}((\partial E\setminus \Gamma)\cap \Omega)=0.$$
\end{itemize}
\end{theorem}
Let us briefly describe the organization of this paper. Section \ref{Preliminary notation and definitions} collects known results, notation and preliminary definitions. Moreover, in this section the equivalence between the constrained problem an the penalized problem is proved. As it always happens when different kind of energies compete with each other, the proof of the regularity is based on the study of the interplay between them. In this case we must compare perimeter and bulk energy (see \cite{AFP}, \cite{Lin}).\\
We point out that the H\"older exponent $\frac{1}{2}$ is critical in this respect for solutions $u$ of either $\eqref{P}$ or $\eqref{P_c}$, in the sense that, whenever $u\in C^{0,\frac{1}{2}}$, under appropriate scaling, the bulk term  locally has the same dimension $n-1$ as the perimeter term.\\
In section 3 we prove suitable energy decay estimates for the bulk energy. The key point of this approach is contained in Lemma \ref{Lemma decadimento 1}, where it is proved that the bulk energy decays faster than $\rho^{n-1}$, that is, for any $\mu\in(0,1)$,
\begin{equation}\label{intro2}
\int_{B_{\rho}(x_0)}|\D u|^2\,dx\leq C \rho^{n-\mu},
\end{equation} 
either in the case that 
$$\min\{|E\cap B_{\rho}(x_0)|,|B_{\rho}(x_0)\setminus E|\}<\varepsilon_0|B_{\rho}(x_0)|,
$$
or in the case that there exists an half-space $H$ such that
$$
|(E\Delta H)\cap B_{\rho}(x_0)|\leq \varepsilon_0|B_{\rho}(x_0)|,
$$
for some $\varepsilon_0>0$. The latter case is the hardest one to handle because it relies on the regularity properties of solutions of a transmission problem which we study in subsection \ref{A decay estimate for elastic minima}. 
Let us notice that, for any given $E\subset \Omega$, local minimizers $u$ of the functional
\begin{equation}\label{intro3}
\int_\Om\bigl [ F(x,u,\nabla u)+\mathbbm{1}_{E}G(x,u,\nabla u)\bigr ]\,dx
\end{equation}
are H\"older continuous, $u\in C^{0,\sigma}_{loc}(\Omega)$, but the needed bound $\sigma>\frac{1}{2}$ cannot be expected in the general case without any information on the set $E$. \\
In subsection \ref{A decay estimate for elastic minima} we prove that minimizers of the functional \eqref{intro3} are in $C^{0,\sigma}$ for every $\sigma\in(0,1)$, in the case $E$ is an half-space. In this context the linearity of the equation strongly comes into play ensuring that the derivatives of the Euler-Lagrange equation are again solutions of the same equation. For the proof in section \ref{Decay of the bulk energy} we readapt a technique depicted in the book \cite{AFP} in the context of the Mumford-Shah functional and recently used in a paper by E. Mukoseeva and G. Vescovo, \cite{MV}.\\
In section \ref{Energy density estimates}, using the estimates obtained in section \ref{Decay of the bulk energy},  we are in position to prove some decay estimates for the whole energy including the perimeter term. More precisely, whenever the  perimeter of $E$ is sufficiently small in a ball $B_{\rho}(x_0)$, then the total energy 
$$
\int_{B_r(x_0)}|\D u|^2\,dx + P(E;B_{r}(x_0)), \quad \quad 0<r<\rho,
$$
decays as $r^{n}$ (see Lemma \ref{Lemma decadimento 2}). 
In the subsequent sections we collect the preliminary results needed to deduce that $\partial E$ is locally represented by a Lipschitz graph, see Theorem \ref{LipApp}.\\
In section \ref{Energy density estimates}, making use of the previous results, we are in position to prove the density upper bound and the density lower bound for the perimeter of $E$ which, in turn, are crucial to prove the Lipschitz approximation theorem. In the subsequent sections the proof strategy follows the path traced from the regularity theory for perimeter minimizers. \\
In section \ref{Compactness for sequences of minimizers} it is proved the compactness for sequences of minimizers which follows in a quite standard way from the density lower bound.\\
Section \ref{Height bound and Lipschitz approximation} is devoted to the Lipschitz approximation theorem which involves the usual main ingredient of the regularity proof, that is the excess
$$
{\mathbf e}(x,r)=\inf_{\nu\in \mathbb{S}^{n-1}}{\mathbf e}(x,r,\nu):= \inf_{\nu\in \mathbb{S}^{n-1}}\frac{1}{r^{n-1}}\int_{\partial E\cap B_r(x)}\frac{|\nu_E(y)-\nu|^2}{2}d\mathcal H^{n-1}(y).
$$ 
In section \ref{Reverse Poincaré inequality} we prove the reverse Poincar\'{e} inequality which is the counterpart of the well-known Caccioppoli's inequality for weak solutions of elliptic equations.\\
Sections \ref{Energy first variation} contains a Taylor-like expansion formula for the terms appearing in the energy under a small domain perturbation.\\
In section \ref{Excess improvement} we finally  prove the excess improvement, which is the main ingredient to achieve the regularity of the interface. More precisely, we prove that,  whenever the excess ${\mathbf e}(x,r)$ goes to zero, for $r\rightarrow 0$, the Dirichlet integral $\int_{B_{\rho}(x_0)}|\D u|^2\,dx$ decays as in \eqref{intro2}. With all these results in hand we can conclude the desired result.\\
In section \ref{Proof of the main theorem} we provide the proof of Theorem \ref{Teorema principale} that is a consequence of the excess improvement proved before.

\section{Preliminary notation and definitions}
\label{Preliminary notation and definitions}
In the rest of the paper we will write $\langle \xi, \eta \rangle$ for the inner product of vectors $\xi, \eta \in \mathbb{R}^n$, and consequently $|\xi|:=\langle \xi, \xi \rangle^{\frac 12}$ will be the corresponding Euclidean norm. As usual $\omega_n$ stands for the Lebesgue measure of the unit ball in $\R^n$.\\
We will denote by $p:\R^n\rightarrow \R^{n-1}$ and $q:\R^n\rightarrow \R$ the horizontal and vertical projections, so that $x=(px,qx)$ for all $x\in \R^n$.
For simplicity of notation we will often write $px=x'$ and $qx=x_n$, so that 
we will write $x=(x',x_n)$, where $x'\in \mathbb{R}^{n-1}$ and $x_n\in \R$.
Accordingly, we denote $\D'=(\partial_{x_1},\dots,\partial_{x_{n-1}})$ the gradient with respect to the first $n-1$ components.\\  
The $n$-dimensional ball in $\R^n$ with center $x_0$ and radius $r>0$ will be denoted as
$$B_R(x_0)=\{x\in\R^{n}: |x-x_0|<R\}.$$ 
If $x_0=0$, we will simply write $B_R$ instead of $B_R(x_0)$.\\ 
The $(n-1)$-dimensional ball in $\R^{n-1}$ with center $x'_0$ and radius $r>0$ will be denoted with a different letter, that is 
$$D_R(x_0)=\{x'\in\R^{n-1}: |x'-x'_0|<R\}.$$
If $u$ is integrable in $B_R(x_0)$ we set
$$
u_R=\frac{1}{\omega_n R^n}\int_{B_R(x_0)} u\,dx = \fint_{B_{R}(x_0)} u\,dx.
$$
For any $\mu\geq 0$ we define the Morrey space $L^{2,\mu}(\Omega)$ as
\begin{equation}
\label{SpazioMorrey}
L^{2,\mu}(\Omega):=\left\{u\in L^2(\Omega)\,:\,\sup_{x_0\in\Omega,\,r>0}r^{-\mu}\int_{\Omega\cap B_r(x_0)}|u|^2\,dx<\infty\right\}.
\end{equation}
In the sequel we will constantly need to denote the difference between $\alpha$ and $\frac{n-1}{n}$, so that we define
\begin{equation}
\label{beta}
\gamma:=n\Bigl(\alpha-\frac{n-1}{n}\Bigr)=1+n(\alpha-1)\in(0,1).
\end{equation}
The following definition is standard.
\begin{definition}
Let $v\in H^1_{loc}(\Omega)$ and assume that $E\subset \Omega$ is fixed. We define the functional $\mathcal{F}_E$ by setting
\begin{equation*}
\mathcal{F}_E(w,\Omega):=\mathcal{F}(E,w;\Omega), \quad\forall w\in H^1(\Omega).
\end{equation*}
Furthermore we say that $v$ is a local minimizer of the integral functional $\mathcal{F}_E$ if and only if
$$
{\mathcal F}_E(v;B_R(x_0))=\min_{w\in v+H^1_0(B_R(x_0))}\mathcal F_E(w;B_R(x_0)),
$$
for all $B_R(x_0)\subset \subset \Omega$.
\end{definition}

It is worth mentioning that for a quadratic integrand $F(x,s,z)$ of the type given in $\eqref{structure1}$ the following growth condition can be immediately deduced from assumptions $\eqref{ellipticity1}$ and $\eqref{ellipticity2}$: 
\begin{equation}\label{ggrowth}
\frac{\nu}{2}|z|^2-\frac{L^2}{\nu}\leq F(x,s,z)\leq (N+1)|z|^2+L(L+1),\quad\forall x\in\Omega,\,\forall s\in\R,\,\forall z\in\R^n.
\end{equation}
The next lemma is very standard and can be found for example in \cite[Lemma 7.54]{AFP}.
\begin{lemma}
\label{Lemma iterativo 2}
Let $f:(0,a]\rightarrow [0,\infty)$ be an increasing function such that
$$
f(\rho)\leq A\Bigl[\Bigl(\frac{\rho}{R}\Bigr)^p+R^s\Bigr]f(R)+BR^q,\quad\text{whenever }0<\rho<R\leq a,
$$
for some constants $A,B\geq 0$, $0<q<p$, $s>0$. Then there exist $R_0=R_0(p,q,s,A)$ and $c=c(p,q,A)$ such that
$$
f(\rho)\leq c\Bigl(\frac{\rho}{R}\Bigr)^qf(R)+ cB\rho^q, \qquad\mbox{whenever }0<\rho<R\leq \min\{R_0,a\}.
$$
\end{lemma}
\subsection{From constrained to penalized problem}
The next theorem allows us to overcome the difficulty of handling with the constraint $|E|=d$. Indeed, we prove that every minimizer of the constrained problem $\eqref{P_c}$ is also a minimizer of a suitable unconstrained problem with a volume penalization of the type given in $\eqref{P}$.
\begin{theorem}
\label{Teorema Penalizzazione} There exists $\Lambda_{0}>0$ such that if $(E,u)$ is a minimizer of the functional
\begin{equation}
\label{Penalized}
{\mathcal F}_{\Lambda}(A,w)=\int_\Om\big [F(x,w,\nabla w)+\mathbbm{1}_{A}G(x,w,\nabla w)\,dx\big]\,dx +P(A;\Om)+\Lambda\big||A| - d\big|^\alpha,
\end{equation}
for some $\Lambda \geq \Lambda_0$, among all configurations $(A,w)$ such that $w=u_0$ on $\partial \Omega$,
then $|E|=d$ and $(E,u)$ is a minimizer of problem $\eqref{P_c}$.
Conversely, if $(E,u)$ is a minimizer of problem $\eqref{P_c}$, then it is a minimizer of \eqref{Penalized}, for all $\Lambda \geq \Lambda_0$.
\end{theorem}
\begin{proof}
The proof can be carried out as in \cite[Theorem 1]{EF}. For reader's convenience we give here its sketch, emphasizing main ideas and minor differences with respect to the case treated in \cite{EF}.\\
The first part of the theorem can be proved by contradiction. Assume that there exist a sequence $(\lambda_h)_{h\in \mathbb N}$ such that $\lambda_h\rightarrow \infty$ as $h\rightarrow \infty$ and a sequence of configurations $(E_h,u_h)$ minimizing $\mathcal{F}_{\lambda_h}$ and such that $u_h=u_0$ on $\partial \Omega$ and $|E_h|\neq d$ for all $h\in\N$. Let us choose now an arbitrary fixed $E_0\subset \Omega$ with finite perimeter such that $|E_0|=d$. Let us point out that 
\begin{equation}\label{Theta}
\mathcal{F}_{\lambda_h}(E_h,u_h)\leq\mathcal{F}(E_0,u_0):=\Theta.
\end{equation}
Without loss of generality we may assume that $|E_h|<d$. Indeed, the case $|E_h|>d$ can be treated in the same way considering the complement of $E_h$ in $\Omega$. Our aim is to show that, for $h$ sufficiently large, there exists a configuration $(\widetilde{E}_h,\tilde{u}_h)$ such that $\mathcal{F}_{\lambda_h}(\widetilde{E}_h,\tilde{u}_h)< \mathcal{F}_{\lambda_h}({E_h},{u_h})$, thus proving the result by contradiction.\\
By condition $\eqref{Theta}$, it follows that the sequence $(u_h)_h$ is bounded in $H^1(\Omega)$, the perimeters of the sets $E_h$ in $\Omega$ are bounded and $|E_h|\rightarrow d$. Therefore, possibly extracting a not relabelled subsequence, we may assume that there exists a configuration $(E,u)$ such that $u_h\rightarrow u$ weakly in $H^1(\Omega)$, $\mathbbm{1}_{E_h}\rightarrow \mathbbm{1}_{E}$ a.e. in $\Omega$, where the set $E$ is of finite perimeter in $\Omega$ and $|E|=d$. The couple $(E,u)$ will be used as reference configuration for the definition of $(\widetilde{E}_h,\tilde{u}_h)$.\\

{\bf Step 1.} {\em Construction of $(\widetilde{E}_h,\tilde{u}_h)$}.
Proceeding exactly as in \cite{EF}, we take a point $x\in\partial^*E\cap\Om$ and observe that  the sets $E_r=(E-x)/r$ converge locally in measure to the half-space $H=\{\langle z,\nu_E(x)\rangle<0\}$, i.e., $\mathbbm{1}_{E_r}\to\mathbbm{1}_H$ in $L^1_{\rm loc}(\R^n)$, where $\nu_E(x)$ is the generalized exterior normal to $E$ at $x$ (see \cite[Definition 3.54]{AFP}). Let $y\in B_1(0)\setminus H$ be the point $y=\nu_E(x)/2$. Given $\e$ (that will be chosen in the Step {2}), since $\mathbbm{1}_{E_r}\to\mathbbm{1}_H$ in $L^1(B_1(0))$ there exists $0<r<1$ such that
$$
|E_r\cap B_{1/2}(y)|<\e,\qquad |E_r\cap B_1(y)|\geq|E_r\cap B_{1/2}(0)|>\frac{\omega_n}{2^{n+2}}\,,
$$
where $\omega_n$ denotes the measure of the unit ball of $\R^n$. Then, if we define $x_r:=x+ry\in\Om$, we have that
$$
|E\cap B_{r/2}(x_r)|<\e r^n,\qquad|E\cap B_r(x_r)|>\frac{\omega_nr^n}{2^{n+2}}\,.
$$
Let us assume, without loss of generality, that $x_r=0$. From the convergence of $E_h$ to $E$ we have that for all $h$ sufficiently large
\begin{equation}\label{unodue}
|E_h\cap B_{r/2}|<\e r^n,\qquad|E_h\cap B_r|>\frac{\omega_nr^n}{2^{n+2}}\,.
\end{equation}
Let us now define the following bi-Lipschitz function used in \cite{EF} which maps $B_r$ into itself: 
\begin{equation}\label{unotre}
\Phi(x)=
\begin{cases}
\bigl(1-\sigma_h(2^n-1)\bigr)x & \text{if}\,\,\,|x|<\displaystyle\frac{r}{2},\\
x+\sigma_h\Bigl(\displaystyle1-\frac{r^n}{|x|^n}\Bigr)x & \text{if}\,\,\,\displaystyle\frac{r}{2}\leq|x|<r,\\
x & \text{if}\,\,\,|x|\geq r\,,
\end{cases}
\end{equation}
for some $0<\sigma_h<1/2^n$ sufficiently small to be chosen later
in such a way that, setting
$$
{\widetilde E}_h:=\Phi(E_h),\qquad{\tilde u}_h:=u_h\circ\Phi^{-1},
$$
we have 
$$|{\widetilde E}_h|<d.$$
We are going to evaluate
\begin{align}\label{unoquattro}
&{\mathcal F}_{\lambda_h}(E_h,u_h)-{\mathcal F}_{\lambda_h}({\widetilde E}_h,{\tilde u}_h)
\nonumber =\biggl[\int_{B_r}\big[F(x,u_h,\nabla u_h)+\mathbbm{1}_{E_h}G(x,u_h,\nabla u_h)\big]\,dx\\
&-\int_{B_r}\big[F(x,\tilde{u}_h,\nabla \tilde{u}_h)+\mathbbm{1}_{\widetilde{E}_h}G(x,\tilde{u}_h,\nabla \tilde{u}_h)\big]\,dy\biggr] 
\\
& \quad+\bigl[P(E_h;{\overline B}_r)-P({\widetilde E}_h;\overline{B}_r)\bigr]+\lambda_h\bigl[(d-|{ E}_h|)^\alpha-(d-|\widetilde{E}_h|)^\alpha\bigr] \nonumber \\
&= I_{1,h}+I_{2,h}+I_{3,h}. \nonumber
\end{align}

{In order to estimate the contribution of the last integrals we need some preliminary estimates for the map $\Phi$ that can be obtained by direct computation (see \cite{EF} or \cite{EL} for the explicit calculation). 
We just observe that for $|x|<r/2$, $\Phi$ is simply a homothety and all the estimates that we are going to introduce are trivial.\\ Conversely, for $r/2<|x|<r$ we have
\begin{equation}\label{JacobianPhi}
\frac{\partial \Phi_i}{\partial x_j}(x)=\Bigl(1+\sigma_h-\frac{\sigma_h r^n}{|x|^n}\Bigr)\delta_{ij}+n\sigma_h r^n\frac{x_ix_j}{|x|^{n+2}}.
\end{equation}
It is clear from this expression that, since $\sigma_h$ is going to zero, $\nabla\Phi$ is a small perturbation of the identity that can be written as 
$$
\nabla \Phi=Id+\sigma_h Z.
$$
We can also address the reader to section $17.2$ ``Taylor's expansion of the determinant close to the identity'' in \cite{Ma} for related estimates.
Then we have 
\begin{equation}\label{Def}
|z-z\circ\nabla \Phi(y)|\leq C_1(n)\sigma_h |z|,\quad
\text{for all}\,\,\, y, z\in \R^n.
\end{equation}
}
It is not difficult to find out also that
\begin{equation}\label{phi}
\bigl\|\nabla\Phi^{-1}\bigl(\Phi(x)\bigr)\bigr\|_\infty\leq\bigl(1-(2^n-1)\sigma_h\bigr)^{-1}\leq1+2^nn\sigma_h, \qquad\text{for all}\,\,\,x\in B_r.
\end{equation}
Concerning $J\Phi$, the Jacobian of $\Phi$, from $\eqref{JacobianPhi}$ we deduce
$$
J\Phi(x)=\Bigl(1+\sigma_h+\frac{(n-1)\sigma_h r^n}{|x|^n}\Bigr)\Bigl(1+\sigma_h-\frac{\sigma_h r^n}{|x|^n}\Bigr)^{{n-1}}.
$$
{For  $r/2<|x|<r$}, we can estimate (see also section 3 in \cite{CFP}):
\begin{align*}
J\Phi(x)
& \geq \Bigl(1+\sigma_h+\frac{(n-1)\sigma_h r^n}{|x|^n}\Bigr)\Bigl(1+\sigma_h-(n-1)\frac{\sigma_h r^n}{|x|^n}\Bigr)\\
& \geq 1+2\sigma_h -\bigl(4^n(n-1)^2-1\bigr)\sigma_h^2>1+\sigma_h,
\end{align*}
provided that we choose  
$$\sigma_h<\frac{1}{4^n(n-1)^2-1}.$$
Summarizing we gain the following inequalities for the Jacobian of $\Phi$:
{
\begin{align}
\label{Jphi}
&1+\sigma_h\leq J\Phi(x),\quad\text{for all }x\in B_r\setminus B_{r/2},\\
& J\Phi(x)\leq1+2^nn\sigma_h,\quad\text{for all }x\in B_r.
\end{align}
Now, let us start estimating $I_{3,h}$ thus proving at the same time that the condition 
$|{\widetilde E}_h|<d$ is satisfied.\\

{\bf Step 2.} {\em Estimate of $I_{3,h}$}.
First we recall \eqref{unodue}, \eqref{unotre}, \eqref{Jphi}, thus getting}
\begin{align*}
|\tilde{E}_h|-|E_h|
& = \int_{E_h\cap B_r\setminus B_{r/2}}\!\left(J\Phi(x)-1\right)\,dx+\int_{E_h\cap B_{r/2}}\!\left(J\Phi(x)-1\right)\,dx\\
& \geq\Bigl(\frac{\omega_n}{2^{n+2}}-\e\Bigr)\sigma_h r^n-\bigl[1-\bigl(1-(2^n-1)\sigma_h\bigr)^n\bigr]\e r^n \\
& \geq \sigma_h r^n\Bigl[\frac{\omega_n}{2^{n+2}}-\e-(2^n-1)n\e\Bigr].
\end{align*}
Therefore, if we choose $0<\e<\e_0(n)$, we have that
\begin{equation}\label{I3hh}
\lambda_h(|\tilde{E}_h|-|E_h|)\geq\lambda_h C_2(n)\sigma_h r^n.
\end{equation}
Moreover, if we denote $\delta_h:=d-|E_h|$, we choose $\sigma_h$ in such a way that
$|\tilde{E}_h|-|E_h|\leq \delta_h/2 $ thus respecting the condition $|{\widetilde E}_h|<d$.
For this reason let us observe that we have, proceding as before and using \eqref{Jphi},
\begin{align*}
|\tilde{E}_h|-|E_h|
& = \int_{E_h\cap B_{r}}\!\left(J\Phi(x)-1\right)\,dx
\leq n2^n \sigma_h r^n.
\end{align*}
Then we will choose 
$$
\delta_h\leq \sigma_h\leq \frac{\delta_h}{n2^{n+1}r^n}.
$$
Let us observe that in the last condition we imposed also that $\sigma_h$ is comparable with $\delta_h$, which is crucial in the following estimate.
Resuming \eqref{I3hh} we can conclude
\begin{align}\label{I3h}
I_{3,h}
& =\lambda_h\bigl[(d-|{ E}_h|)^\alpha-(d-|\widetilde{E}_h|)^\alpha\bigr]\geq \lambda_h\frac{\alpha}{(d-|E_h|)^{1-\alpha}}(|\tilde{E}_h|-|E_h|)\\
& = \lambda_h\alpha(d-|E_h|)^{\alpha}\frac{|\tilde{E}_h|-|E_h|}{d-|E_h|}\geq \lambda_h\alpha \delta_h^\alpha\frac{C_2(n)\sigma_h r^n}{\delta_h}\\
& \geq \lambda_h C_3(n,\alpha)\sigma_h^\alpha r^n,
\end{align}
for some positive constant $C_3=C_3(n,\alpha)$.\\
{\bf Step 3.} {\em Estimate of $I_{1,h}$}. 
Now we can perform the change of variables $y=\Phi(x)$ and, observing that $\mathbbm{1}_{\widetilde{E}_h}(\Phi(x))=\mathbbm{1}_{E_h}(x)$, we get
\begin{align*}
I_{1,h}
& =\int_{B_r}\bigl[F(x,u_h,\nabla u_h) - J\Phi(x) F(\Phi(x),u_h(x),\nabla u_h(x)\circ \nabla \Phi^{-1}(\Phi(x)))\bigr]\,dx\\  
& +\int_{B_r\cap E_h}\bigl[G(x,u_h,\nabla u_h) - J\Phi(x) G(\Phi(x),u_h(x),\nabla u_h(x)\circ \nabla \Phi^{-1}(\Phi(x)))\bigr]\,dx\\
& =:J_{1,h}+J_{2,h}. 
\end{align*}
The two terms $J_{1,h}$ and $J_{2,h}$, involving $F$ and $G$ in $B_r$ and $B_r\cap E_h$ respectively, can be treated in the same way. Therefore we just perform the calculation for $J_{1,h}$.\\ To make the argument more clear, since we shall use the structure conditions $\eqref{structure1}$ and $\eqref{structure2}$ we introduce the following notation. $A_2(x,s)$ denotes the quadratic form and $A_1(x,s)$ denotes the linear form defined as follows:
\begin{align*}
& A_2(x,s)[z]:=a_{ij}(x,s)z_i z_j,\quad A_1(x,s)[z]:=a_{i}(x,s)z_i,
\end{align*}
for any $z\in\R^n$. Analogously we set $A_0(x,s)=a(x,s)$. Accordingly, we can write down
\begin{align}\label{J1h}
& J_{1,h}\nonumber\\ 
& =\int_{B_r} \!\Bigl\{A_2(x,u_h(x))[\nabla u_h(x)]\!-\!A_2(\Phi(x),u_h(x))[\nabla u_h(x)\!\circ\! \nabla \Phi^{-1}(\Phi(x))]J\Phi(x)\Big\}\,dx \nonumber\\
& +\int_{B_r} \!\Bigl\{A_1(x,u_h(x))[\nabla u_h(x)]\!-\!A_1(\Phi(x),u_h(x))[\nabla u_h(x)\!\circ\! \nabla \Phi^{-1}(\Phi(x))]J\Phi(x)\Big\}\,dx \nonumber\\
& +\int_{B_r}\!\Bigl\{A_0(x,u_h(x))\!-\!A_0(\Phi(x),u_h(x))J\Phi(x)\Big\}\,dx.
\end{align}
We proceed estimating the first difference in the previous equality, being the other similar and indeed easier to handle.
\begin{align*}
& \int_{B_r} \!\Bigl\{A_2(x,u_h(x))[\nabla u_h(x)]\!-\!A_2(\Phi(x),u_h(x))[\nabla u_h(x)\!\circ\! \nabla \Phi^{-1}(\Phi(x))]J\Phi(x)\Big\}\,dx \\
& =\int_{B_r} \!\Bigl\{A_2(\Phi(x),u_h(x))[\nabla u_h(x)]-A_2(\Phi(x),u_h(x))[\nabla u_h(x)\!\circ\! \nabla \Phi^{-1}(\Phi(x))]J\Phi(x)\Big\}\,dx\\
&+\int_{B_r} \!\Bigl\{A_2(x,u_h(x))[\nabla u_h(x)]-A_2(\Phi(x),u_h(x))[\nabla u_h(x)]\Big\}\,dx =: H_{1,h}+ H_{2,h}.
\end{align*}
The first term $H_{1,h}$ can be estimated observing that, as a consequence of  $\eqref{ellipticity1}$, we have:
$$
|A_2[\xi]-A_2[\eta]|\leq N|\xi+\eta||\xi-\eta|,\quad \forall\xi,\eta\in \mathbb{R}^n.
$$
If we apply the last inequality to the vectors
$$
\xi:=\nabla u_h(x), \qquad \qquad \eta:={\sqrt{J\Phi(x)}}[\nabla u_h(x)\!\circ\! \nabla \Phi^{-1}(\Phi(x))],
$$
we are led to estimate $|\xi-\eta|$.\\
We start observing that, being $J\Phi(x)=\bigl(1-\sigma_h(2^n-1)\bigr)^n$ for $|x|<r/2$, by also using, $\eqref{Jphi}$ we deduce
$$
|\sqrt{J\Phi(x)}-1|<C(n)\sigma_h,\quad \text{ for all }x\in \R^n.
$$
Therefore we have
$$
{\sqrt{J\Phi}} \xi-\xi|\leq C(n)\sigma_h|\xi|.
$$
In addition  choosing $z=\xi\circ \nabla\Phi^{-1}(\Phi(x))$ in $\eqref{Def}$ and using also $\eqref{phi}$, we can deduce 
$$
|\xi \circ \nabla \Phi^{-1}(\Phi(x))-\xi|\leq \sigma_h C_1(n)|\xi\circ \nabla \Phi^{-1}(\Phi(x))|\leq \sigma_h |\xi|C_1(n) \norm{\nabla \Phi^{-1}(\Phi(x))}_\infty\leq n2^nC_1(n)\sigma_h|\xi|.
$$
Summarizing we finally get
$$
|\xi-\eta|\leq \sigma_h C(n)|\nabla u_h(x)|,\qquad \qquad|\xi+\eta|\leq C(n)|\nabla u_h(x)|,
$$
for some constant $C=C(n)>0$. From the previous estimates we deduce that
\begin{equation}\label{H1h}
|H_{1,h}|\leq \sigma_h N C^2(n) \int_{B_r}|\nabla u_h(x)|^2\,dx\leq  \sigma_h N C^2(n)\Theta,
\end{equation}
where $\Theta$ is defined in \eqref{Theta}. The second term $H_{2,h}$ can be estimated using the H\"older continuity assumption on $a_{ij}$ and observing that $|x-\Phi(x)|\leq \sigma_h r2^n$. Therefore we deduce that
\begin{equation}\label{H2h}
|H_{2,h}|\leq (\sigma_h r2^n)^\alpha L_\alpha \int_{B_r}|\nabla u_h(x)|^2\,dx\leq  \sigma_h^\alpha  C(n,\alpha,L_\alpha)\Theta.
\end{equation}
In conclusion, since the other terms in $\eqref{J1h}$ can be estimated in the same way,
collecting estimates $\eqref{H1h}$ and $\eqref{H2h}$ we get
$$
|J_{1,h}|\leq \sigma_h^\alpha  C(n,N,\alpha,L_\alpha)\Theta.
$$
Since the same estimate holds true for $J_{2,h}$, we conclude that
\begin{equation}
\label{I1h}
I_{1,h}\geq -\sigma_h^\alpha C_4(n,N,\alpha,L_\alpha)\Theta,
\end{equation}
for some constant $C_4=C_4(n,N,\alpha,L_\alpha)>0$.\\
{\bf Step 4.} {\em Estimate of $I_{2,h}$}. In order to estimate $I_{2,h}$, we can use the area formula for maps between rectifiable sets. If we denote by $T_{h,x}$ the tangential gradient of $\Phi$ along the approximate tangent space to $\partial^* E_h$ in $x$ and $T^*_{h,x}$ is the adjoint of the map $T_{h,x}$, the $(n-1)$-dimensional jacobian of $T_{h,x}$ is given by
$$
J_{n-1}T_{h,x}=\sqrt{{\rm det}\bigl(T^*_{h,x}\circ T_{h,x}\bigr)}.
$$
Thereafter we can estimate
\begin{equation}\label{TJ}
J_{n-1}T_{h,x}\leq 1+\sigma_h+2^n(n-1)\sigma_h.
\end{equation}
We address the reader to \cite{EF} where explicit calculations are given. In order to estimate $I_{2,h}$, we use the area formula for maps between rectifiable sets (\cite[Theorem~2.91]{AFP}), thus getting
\begin{align*}
I_{2,h}
& =P(E_h;{\overline B}_r)-P({\widetilde E}_h;{\overline B}_r)=\int_{\partial^*E_h\cap{\overline B}_r}\!d\H^{n-1}-\int_{\partial^*E_h\cap{\overline B}_r}\!J_{n-1}T_{h,x}\,d\H^{n-1} \\
& =\int_{\partial^*E_h\cap{\overline B}_r\setminus B_{r/2}}\!\left(1-J_{n-1}T_{h,x}\right)\,d\H^{n-1}+\int_{\partial^*E_h\cap B_{r/2}}\!\left(1-J_{n-1}T_{h,x}\right)\,d\H^{n-1}\,.
\end{align*}
Notice that the last integral in the above formula is non-negative since $\Phi$ is a contraction in $B_{r/2}$, hence $J_{n-1}T_{h,x}<1$ in $B_{r/2}$, while from \eqref{TJ} we have
$$
\int_{\partial^*E_h\cap{\overline B}_r\setminus B_{r/2}}\!\left(1-J_{n-1}T_{h,x}\right)\,d\H^{n-1}\geq-2^nnP(E_h;{\overline B}_r)\sigma_h\geq-2^nn\Theta\sigma_h^\alpha\,,
$$
thus concluding that
\begin{equation}\label{I2h}
I_{2,h}\geq-2^nn\Theta\sigma_h^\alpha.
\end{equation}
{Finally to conclude the proof we recall \eqref{unoquattro}, \eqref{I3h}, \eqref{I1h} and \eqref{I2h}  to obtain}
$$
{\mathcal F}_{\lambda_h}(E_h,u_h)-{\mathcal F}_{\lambda_h}({\widetilde E}_h,{\tilde u}_h)\geq\sigma_h^\alpha\bigl(\lambda_h C_3(n,\alpha)r^n-\Theta(C_4(n,N,\alpha,L_\alpha)+2^nn)\bigr)>0,
$$
if $\lambda_h$ is sufficiently large. This contradicts the minimality of $(E_h,u_h)$, thus concluding the proof.

\end{proof}

The previous theorem  motivates the following definition.
\begin{definition}[$(\Lambda,\alpha)$-minimizers]
The energy pair $(E,u)$ is a $(\Lambda,\alpha)$-minimizer in $\Omega$ of the functional $\mathcal {F}$, defined in \eqref{intro0}, if and only if for every $B_r(x_0)\subset \Omega$ it holds:
$$
\mathcal{F}(E,u;B_r(x_0))\leq\mathcal{F}(F,v;B_r(x_0))+\Lambda |F\Delta E|^\alpha,
$$
whenever $(F,v)$ is an admissible test pair, namely, $F$ is a set of finite perimeter with $F\Delta E\subset \subset B_r(x_0)$ and $v-u\in H^1_0(B_r(x_0))$.
\end{definition}

\section{Decay of the bulk energy}
\label{Decay of the bulk energy}
We start quoting higher integrability results both for local minimizers of the functional $\eqref{intro0}$  and  for comparison functions that we will use later in the paper. We assume that $E$ is fixed and therefore we consider only the dependence on the bulk term through $u$.
It is worth mentioning that the following lemmata can be applied in general to minimizers of integral functionals of the type
\begin{equation}\label{H}
\mathcal{H}(u;\Omega):=\int_{\Omega}H(x,u,\D u)\,dx,
\end{equation}
assuming that the energy density $H$ satisfies only the structure condition $\eqref{structure1}$ and the growth conditions $\eqref{ellipticity1}$ and $\eqref{ellipticity2}$, without assuming any continuity on the coefficients. It is clear that functionals of the type $\eqref{intro0}$ belong to this class and in addition the involved estimates only depend on the constants appearing in $\eqref{ellipticity1}$ and $\eqref{ellipticity2}$ but do not depend on $E$ accordingly. Since the argument is very standard we address the reader to \cite{EL} where detalied proofs is given.

\begin{lemma}
\label{Lemma maggiore sommabilità}
Let $u\in H^1(\Omega)$ be a local minimizer of the functional $\mathcal{H}$ defined in \eqref{H}, where $H$ satisfies the structure condition $\eqref{structure1}$ and the growth conditions $\eqref{ellipticity1}$ and $\eqref{ellipticity2}$. There exists $s=s(n,\nu,N,L)>1$ such that, for every $B_{2R}(x_0)\subset \subset \Omega$, it holds
\begin{equation*}
\fint_{B_{R}(x_0)}|\nabla u|^{2s}\,dx\leq C_1\biggl( \fint_{B_{2R}(x_0)}\big(1+|\nabla u|^{2}\big)\,dx\biggr)^s,
\end{equation*}
where $C_1=C_1(n,\nu,N,L)$ is a positive constant.
\end{lemma}
In the next subsection we will prove some energy  density estimates by using a standard comparison argument. For this purpose we will need a reverse H\"{o}lder inequality for the comparison function defined below.
\begin{definition}[Comparison function]
Let $u\in H^1(\Omega)$ be a local minimizer of the functional $\mathcal{H}$ defined in \eqref{H} and $B_{2R}(x_0)\subset \subset \Omega$. We shall denote by $v$ the solution of the following problem
\begin{equation}\label{Comparison}
v:=\operatornamewithlimits{argmin}_{w\in u+H^1_0(B_R(x_0))}
\int_{B_{R}(x_0)}\tilde{H}(x,\nabla w)\,dx,
\end{equation}
where $\tilde{H}(x,z):=H(x,u(x),z)$ satisfies the structure condition 
$\eqref{structure1}$ and the growth conditions $\eqref{ellipticity1}$ and $\eqref{ellipticity2}$.
\end{definition}
\begin{lemma}
\label{Lemma maggiore sommabilità funzione confronto}
Let $u\in H^1(\Omega)$ be a local minimizer of the functional $\mathcal{H}$ defined in \eqref{H}, where $H$ satisfies the structure condition $\eqref{structure1}$ and the growth conditions $\eqref{ellipticity1}$ and $\eqref{ellipticity2}$. Let $v\in H^{1}(B_R(x_0))$ be the comparison function defined in $\eqref{Comparison}$. Denoting by $s=s(n,\nu,N,L)>1$ the same exponent given in Lemma \ref{Lemma maggiore sommabilità}, it holds
\begin{equation*}
\fint_{B_{R}(x_0)}|\nabla v|^{2s}\,dx\leq C_2\biggl( \fint_{B_{2R}(x_0)}\big(1+|\nabla u|^{2}\big)\,dx\biggr)^s,
\end{equation*}
where $C_2=C_2(n,\nu,N,L)$ is a positive constant.
\end{lemma}

\subsection{A decay estimate for elastic minima}
\label{A decay estimate for elastic minima}
In this section we prove a decay estimate for elastic minima that will be crucial for the proof strategy. 
Indeed, we show that if $(E,u)$ is a $(\Lambda,\alpha)$-minimizer of the functional $\mathcal{F}$ defined in $\eqref{intro0}$ and $x_0$ is a point in $\Omega$, where either the density of $E$ is close to $0$ or $1$, or the set $E$ is asymptotically close to a hyperplane, then for $\rho$ sufficiently small we have
$$
\int_{B_{\rho}(x_0)}|\nabla u_E|^2\,dx\leq C\rho^{n-\mu},
$$
for any $\mu\in(0,1)$. 
A preliminary result we want to mention, which will be used later, provides an upper bound for $\mathcal{F}$. The proof is rather standard and is related to the threshold H\"older exponent $\frac 12$ of the function $u$, when $(E,u)$ is either a solution of the constrained problem \eqref{P_c} or a solution of the penalized problem \eqref{P} defined in Section \ref{Introduction and statements}. 
For the proof we address the reader to \cite[Lemma 2.3]{LK} and \cite{FJ}. 
A detailed proof in the case of costrained problems and for functionals satisfying general $p$-polinomial growth is contained in \cite{CFP}.

\begin{theorem}
\label{Teorema upper bound}\label{Energy upper bound}
Let $(E,u)$ be a $(\Lambda,\alpha)$-minimizer of $\mathcal{F}$ in $\Omega$. Then for every open set $U\subset\subset \Omega$ there exists a constant $C_3=C_3\big(n,\alpha,\Lambda,U,\norm{\D u}_{L^2(\Omega)}\big)>0$ such that for every $B_r(x_0)\subset U$ it holds
\begin{equation*}
\mathcal{F}(E,u;B_r(x_0))\leq C_3 r^{n-1}.
\end{equation*}
\end{theorem}

\begin{proof}
Fixing $B_r(x_0)\subset U\subset\subset\Omega$, we compare $(E,u)$ with $(E\setminus B_r(x_0),u)$ thus obtaining
\begin{align*}
\mathcal{F}(E,u;\Omega)
&\leq \mathcal{F}(E\setminus B_r(x_0),u;\Omega)+\Lambda|E\Delta(E\setminus B_r(x_0))\cap \Omega|^\alpha\\
& \leq\mathcal{F}(E\setminus B_r(x_0),u;\Omega)+\Lambda|B_r(x_0)|^\alpha.
\end{align*}
Making $\mathcal{F}$ explicit and getting rid of the common terms, we obtain:
\begin{align}
\label{a36}
\int_{B_r(x_0)\cap E}G(x,u,\D u)\,dx+P(E;B_r(x_0))
& \leq P(E\cap\dd B_r(x_0);\Omega)+c(n,\alpha,\Lambda)r^{n\alpha}\notag\\
& \leq \mathcal{H}^{n-1}(\dd B_r(x_0))+c(n,\alpha,\Lambda)r^{n-1}\notag\\
&\leq c(n,\alpha,\Lambda)r^{n-1}.
\end{align}
Now we want to prove that there exist $\tau\in\big(0,\frac{1}{2}\big)$ and $\delta\in(0,1)$ such that for every $M>0$ there exists $h_0\in\N$ such that, for any $B_r(x_0)\subset U$, we have
\begin{equation*}
\int_{B_r(x_0)}|\D u|^2\leq h_0r^{n-1} \quad\text{or}\quad \int_{B_{\tau r}(x_0)}|\D u|^2\,dx\leq M\tau^{n-\delta}\int_{B_r(x_0)}|\D u|^2\,dx.
\end{equation*}
\textbf{Step 1: }Arguing by contradiction, for $\tau\in\big(0,\frac{1}{2}\big)$ and $\delta\in(0,1)$, we choose $M>\tau^{\delta-n}$ and we assume that, for every $h\in\N$, there exists a ball $B_{r_h}(x_h)\subset U$ such that
\begin{equation}
\label{a34}
\int_{B_{r_h}(x_h)}|\D u|^2\,dx>hr_h^{n-1}
\end{equation}
and
\begin{equation}
\label{a35}
\int_{B_{\tau r_h}(x_h)}|\D u|^2\,dx> M\tau^{n-\delta}\int_{B_{r_h}(x_h)}|\D u|^2\,dx.
\end{equation}
Note that estimates \eqref{a36} and \eqref{a34} yield
\begin{equation}
\label{a44}
\int_{B_{r_h}(x_h)\cap E}|\D u|^2\,dx+P(E;B_{r_h}(x_h))\leq c_0r_h^{n-1}<\frac{c_0}{h}\int_{B_{r_h}(x_h)}|\D u|^2\,dx,
\end{equation}
and so
\begin{equation}
\label{a41}
\int_{B_{r_h}(x_h)\cap E}|\D u|^2\,dx<\frac{c_0}{h}\int_{B_{r_h}(x_h)}|\D u|^2\,dx,
\end{equation}
for some positive constant $c_0$.\\
 \textbf{Step 2: }We will prove our aim by means of a blow-up argument. We set
\begin{equation*}
\varsigma_h^2:=\dashint_{B_{r_h}(x_h)}|\D u|^2\,dx
\end{equation*}
and, for $y\in B_1$, we introduce the sequence of rescaled functions defined as
\begin{equation*}
v_h(y):=\frac{u(x_h+r_hy)-a_h}{\varsigma_h r_h},\quad\text{with}\quad a_h:=\dashint_{B_{r_h}(x_h)}u\,dx.
\end{equation*}
We have ${\D u(x_h+r_hy)}=\varsigma_h\D v_h(y)$ and a change of variable yields
\begin{equation*}
\int_{B_1}|\D v_h(y)|^2\,dy=\frac{1}{\varsigma_h^2}\dashint_{B_{r_h}(x_h)}|\D u(x)|^2\,dx=1.
\end{equation*}
Therefore, there exist a (not relabeled) subsequence of $v_h$ and $v\in H^1(B_1)$ such that $v_h\rightharpoonup v$ in $H^1(B_1)$ and $v_h\rightarrow v$ in $L^2(B_1)$.
Moreover, the semicontinuity of the norm implies
\begin{equation}
\label{a39}
\dashint_{B_1}|\D v(y)|^2\,dy\leq\liminf_{h\rightarrow\infty}\dashint_{B_1}|\D v_h(y)|^2\,dy=1.
\end{equation}
We rewrite the inequalities \eqref{a34}, \eqref{a35} and \eqref{a41}. They become, respectively,
\begin{equation}
\label{a37}
\varsigma_h^2>\frac{h}{r_h},
\end{equation}
\begin{equation}
\label{a38}
\dashint_{B_\tau}|\D v_h(y)|^2\,dy>M\tau^{-\delta},
\end{equation}
\begin{equation}
\label{a42}
\int_{B_1\cap E^*_h}|\D v_h(y)|^2\,dy<\frac{c_0}{h}\int_{B_1}|\D v_h(y)|^2\,dy=\frac{c_0\omega_n}{h}.
\end{equation}
Of course, \eqref{a37} implies that $\varsigma_h\rightarrow\infty$, as $h\rightarrow\infty$.
\textbf{Step 3: }We claim that the $L^2$-norm of $v_h$ converges to the $L^2$-norm of $v$. Consider the sets
\begin{equation*}
E^*_h:=\frac{E-x_h}{r_h}\cap B_1.
\end{equation*}
Since $r_h^{n-1}P(E^*_h;B_1)=P(E;B_{r_h}(x_h))$, by \eqref{a44}, we have that the sequence $\{P(E^*_h;B_1)\}_{h\in\N}$ is bounded. Therefore up a not relabeled subsequence, $\mathbbm{1}_{E_h}\rightarrow\mathbbm{1}_{E^*}$ in $L^1(B_1)$, for some set $E^*\subset B_1$ of locally finite perimeter. By \eqref{a42} and Fatou's Lemma,
\begin{equation}
\label{a43}
\int_{B_1\cap E^*}|\D v(y)|^2\,dy=0.
\end{equation}
By $\Lambda$-minimality of $(E,u)$ with respect to $(E,u+\phi)$ we get, for $\phi\in H^1_0(B_{r_h}(x_h))$,
\begin{align*}
& \int_{B_{r_h}(x_h)}\big[F(x,u,\D u)+\mathbbm{1}_E G(x,u,\D u)\big]\,dx\\
& \leq\int_{B_{r_h}(x_h)}\big[F(x,u+\phi,\D u+\D\phi)+\mathbbm{1}_E G(x,u+\phi,\D u+\D\phi)\big]\,dx.
\end{align*}
Using the change of variable $x=x_h+r_hy$, we deduce for every $\psi\in H^1_0(B_1)$,
\begin{align*}
& \int_{B_1}\big[F(x_h+r_hy,u(x_h+r_hy),\varsigma_h\D v_h)+\mathbbm{1}_{E^*_h} G(x_h+r_hy,u(x_h+r_hy),\varsigma_h\D v_h)\big]\,dy\\
& \leq \int_{B_1}F(x_h+r_hy,u(x_h+r_hy)+r_h\psi,\varsigma_h\D v_h+\D\psi)\,dy\\
&+\int_{B_1}\mathbbm{1}_{E^*_h} G(x_h+r_hy,u(x_h+r_hy)+r_h\psi,\varsigma_h\D v_h+\D\psi)\,dy.
\end{align*}
Let $\eta\in C^\infty_c(B_1)$ such that $0\leq\eta\leq 1$. We choose the test function $\psi_h=\varsigma_h\eta(v-v_h)$ and exploit $\D v_h+\D\psi_h$ for reader convenience,
$$\D v_h+\D\psi_h=\varsigma_h\eta\D v+\varsigma_h(1-\eta)\D v_h+\varsigma_h(v-v_h)\D\eta.$$ 
For simplicity of notation we will denote $w_h:=u(x_h+r_hy)+r_h\varsigma_h\eta(v-v_h)$ so that the previous inequality can be read as
\begin{align*}
& \int_{B_1}\big[F(x_h+r_hy,u(x_h+r_hy),\varsigma_h\D v_h)+\mathbbm{1}_{E^*_h} G(x_h+r_hy,u(x_h+r_hy),\varsigma_h\D v_h)\big]\,dy\\
& \leq\int_{B_1}F(x_h+r_hy,w_h,{{\varsigma_h\eta\D v}}+\varsigma_h(1-\eta)\D v_h+{\color{red}\varsigma_h(v-v_h)\D\eta})\,dy\\
& + \int_{B_1}\mathbbm{1}_{E^*_h}G(x_h+r_hy,w_h,{\varsigma_h\eta\D v}+\varsigma_h(1-\eta)\D v_h+{\color{red}\varsigma_h(v-v_h)\D\eta})\,dy.
\end{align*}
Using the quadratic structure of $F$ and $G$ we can pull out the terms in red in order to use the convexity in the next step.
\begin{align*}
& \int_{B_1}\big[F(x_h+r_hy,u(x_h+r_hy),\varsigma_h\D v_h)+\mathbbm{1}_{E^*_h} G(x_h+r_hy,u(x_h+r_hy),\varsigma_h\D v_h)\big]\,dy\\
& \leq\int_{B_1}F(x_h+r_hy,w_h,\varsigma_h\eta\D v+\varsigma_h(1-\eta)\D v_h)\,dy + \int_{B_1}\mathbbm{1}_{E^*_h}G(x_h+r_hy,w_h,\varsigma_h\eta\D v+\varsigma_h(1-\eta)\D v_h)\,dy\\
& +c(N,L)\int_{B_1}\big(|\varsigma_h\D v|+|\varsigma_h\D v_h|+|\varsigma_h(v-v_h))|\big)\varsigma_h|v-v_h|\,dy.
\end{align*}
Using the convexity of $F$ and $G$ and rearranging the terms we obtain
\begin{align}
\label{a40}
& \int_{B_1}\eta F(x_h+r_hy,w_h,\varsigma_h\D v_h)\leq \int_{B_1} \eta F(x_h+r_hy,w_h,\varsigma_h\D v)\,dy\notag\\
& +\int_{B_1}\big[F(x_h+r_hy,w_h,\varsigma_h\D v_h)-F(x_h+r_hy,u(x_h+r_hy),\varsigma_h\D v_h)\big]\,dy\notag\\
& +\int_{B_1}\mathbbm{1}_{E^*_h}\big[G(x_h+r_hy,w_h,\varsigma_h\D v_h)-G(x_h+r_hy,u(x_h+r_hy),\varsigma_h\D v_h)\big]\,dy\notag\\
& +\int_{B_1}\mathbbm{1}_{E_h^*}\eta\big[G(x_h+r_hy,w_h,\varsigma_h\D v)-G(x_h+r_hy,w_h,\varsigma_h\D v_h)\big]\,dy\\
& +c(N,L)\int_{B_1}\big(|\varsigma_h\D v|+|\varsigma_h\D v_h|+|\varsigma_h(v-v_h))|\big)\varsigma_h|v-v_h|\,dy.
\end{align}
The last term and the second to last term can be treated in a standard way using \eqref{a39},
H\"older's inequality, the strong convergence of $v_h$ to $v$ and the weak convergence of $\D v_h$ to $\D v$. The remaining two terms, which differ only in the second argument, can be treated as follows.\\
We remark that by definition of $v_h$ and H\"older continuity of $u_h$ immediately follows $r_h\varsigma_h v_h\rightarrow 0$. Therefore, being $r_h\varsigma_h\rightarrow 0$ where $v\neq 0$, we deduce also $w_h-u(x_h+r_hy)=r_h\varsigma_h\eta(v-v_h)\rightarrow 0$ for a.e. $y\in B_1$. Finally, using the equi-integrability of $|\nabla v_h|^2$, resulting from the weak convergence of $\nabla v_h$, and the boundedness of the coefficients $a_{ij},a_i,a$ we conclude that
\begin{align*}
& \int_{B_1}\big[F(x_h+r_hy,w_h,\varsigma_h\D v_h)-F(x_h+r_hy,u(x_h+r_hy),\varsigma_h\D v_h)\big]\,dy\\
& \leq \varsigma_h^2\int_{B_1}|a_{ij}(x_h+r_hy,w_h)-a_{ij}(x_h+r_hy,u(x_h+r_hy)|\D_iv_h||\D_jv_h|\,dy\\
& +\varsigma_h\int_{B_1}|a_i(x_h+r_hy,w_h)-a_i(x_h+r_hy,u(x_h+r_hy)|\D_iv_h|\,dy +c(n,L)=\varsigma_h^2\varepsilon_h.
\end{align*}
Combining the previous inequalities, we get
\begin{equation*}
\int_{B_1}\eta F(x_h+r_hy,w_h,\varsigma_h\D v_h)\,dy\leq \int_{B_1}\eta F(x_h+r_hy,w_h,\varsigma_h\D v)\,dy+\varsigma_h^2\varepsilon_h.
\end{equation*}
Dividing by $\varsigma_h^2$, the linear terms in $F$ tend to 0, thus getting
\begin{equation*}
\int_{B_1}\eta a_{ij}(x_h+r_h y,w_h)\D_i v_h\D_j v_h\,dy\leq\int_{B_1}\eta a_{ij}(x_h+r_h y,w_h)\D_i v\D_j v\,dy+\varepsilon_h.
\end{equation*}
Since $B_{r_h}(x_h)\subset U\subset\subset\Omega$ for all $h\in\N$, we may assume that $x_h\rightarrow \overline{x}$, as $h\rightarrow\infty$. Letting $\eta\downarrow 1$ in the previous inequality, passing to the lower limit, as $h\rightarrow\infty$,  by lower semicontinuity, we finally get
\begin{equation*}
\lim_{h\rightarrow\infty}\int_{B_1} a_{ij}(\overline{x},u(\overline{x}))\D_i v_h\D_j v_h\,dy=\int_{B_1} a_{ij}(\overline{x},u(\overline{x}))\D_i v\D_j v\,dy.
\end{equation*}
Since the matrix $a_{ij}(\overline{x},u(\overline{x}))$ is elliptic and bounded, it induces a norm which is equivalent to the euclidean norm. Thus we get
\begin{equation*}
\lim_{h\rightarrow\infty}\dashint_{B_\tau}|\D v_h|^2\,dy=\dashint_{B_\tau}|\D v|^2\,dy\leq\frac{1}{\tau^n}\dashint_{B_1}|\D v|^2\,dy\leq\frac{1}{\tau^n},
\end{equation*}
which contradicts \eqref{a38}, provided we choose $M>\tau^{\delta-n}$.\\
\textbf{Step 4: } We conclude that there exists $\tau\in\big(0,\frac{1}{2}\big)$ and $\delta\in(0,1)$ such that, setting $M=1$, there exists $h_0\in\N$ such that, for any $B_r(x_0)\subset\Omega$, we have
\begin{equation*}
\int_{B_r(x_0)}|\D u|^2\leq h_0r^{n-1} \quad\text{or}\quad \int_{B_{\tau r}(x_0)}|\D u|^2\,dx\leq \tau^{n-\delta}\int_{B_r(x_0)}|\D u|^2\,dx.
\end{equation*}
Hence, 
\begin{equation*}
\int_{B_{\tau r}(x_0)}|\D u|^2\,dx\leq \tau^{n-\delta}\int_{B_r(x_0)}|\D u|^2\,dx+h_0r^{n-1},
\end{equation*}
and, using Lemma \ref{Lemma iterativo 2}, we obtain that
\begin{equation*}
\int_{B_\rho(x_0)}|\D u|^2\,dx\leq c\Bigg\{\bigg(\frac{\rho}{r}\bigg)^{n-1}\int_{B_r(x_0)}|\D u|^2\,dx+h_0\rho^{n-1} \Bigg\},\quad\forall \,0<\rho<r\leq R,
\end{equation*}
and so
\begin{equation*}
\int_{B_\rho(x_0)}|\D u|^2\,dx\leq c\rho^{n-1}.
\end{equation*}
\end{proof}
As a consequence of the previous theorem, using Poincaré's inequality and the characterization of Campanato spaces (see for example \cite[Theorem 2.9]{Giu}), we can infer that $u\in C^{0,\frac 12}$. We deduce the following remark.
\begin{remark}
\label{Osservazione Holderianita}
Let $(E,u)$ be a $\Lambda$-minimizer of the functional ${\mathcal F}$ defined in \eqref{intro0}. For every open set $U\subset \subset \Omega$ there exists a constant $C=C\big(n,\alpha,\Lambda,U,\norm{\D u}_{L^2(\Omega)}\big)>0$ such that 
\begin{equation}\label{hoelderu}
\sup_{x,y \in U}\frac{|u(x)-u(y)|}{|x-y|^{\frac 12}}\leq C.
\end{equation}
\end{remark}

In order to prove the main lemma of this section we introduce the following preliminary result. For reader's convenience we give here a sketch of the proof, which can be found in \cite{MV}. Actually we state here a weaker version that is suitable for our aim. In the following we will denote 
$$H=\lbrace x\in\R\,:\, x_n>0 \rbrace.
$$
\begin{lemma}
\label{Lemma semispazio}
Let $v\in H^1(B_1)$ be a solution of
\begin{equation*}
-\textnormal{div}(A\D u)=\textnormal{div}\,G, \qquad \mbox{ in }\mathcal{D}'(B_1),
\end{equation*}
where
$$
G^+:=\mathbbm{1}_H G\in C^{0,\sigma}(H\cap B_1),\qquad 
G^-:=\mathbbm{1}_{H^c}G\in C^{0,\sigma}(H^c\cap B_1),
$$
for some $\sigma\in(0,1]$ and $A$ is an elliptic matrix satisfying 
$$\nu|z|^2\leq A_{ij}(x)z_i z_j\leq N |z|^2$$ and 
$$
A^+:=\mathbbm{1}_H A\in C^{0,\sigma}(\overline{H}\cap B_1),\qquad 
A^-:=\mathbbm{1}_{H^c}A\in C^{0,\sigma}(\overline{H}^c\cap B_1),
$$
for some constants $\nu,N>0$. Let us denote
$$
C_A=\max\big\{\norm{A^+}_{C^{0,\sigma}},\norm{A^-}_{C^{0,\sigma}}\big\},\qquad C_G=\max\big\{\norm{G^+}_{C^{0,\sigma}},\norm{G^-}_{C^{0,\sigma}}\big\}.
$$
Then $\D v\in L_{loc}^{2,n}(B_1)$ (see \eqref{SpazioMorrey}). Moreover, there exist two constants $C=C\big(n,\nu,N,C_A,C_G\big)$ and $r_0=r_0(n,\nu,N,\norm{G}_{L^\infty},C_A,C_G)$ such that, for any $r<r_0$ with $B_r(x_0)\subset B_1$,
\begin{equation}\label{MVdecay}
\int_{B_{\rho}(x_0)}|\D v|^2\,dx \leq C\Bigl(\frac{\rho}{r}\Bigr)^{n}\int_{B_{r}(x_0)}|\D v|^2\,dx+C\rho^{n}, \quad \forall\, \rho<\frac{r}{4}.
\end{equation}
\end{lemma}
\begin{proof}
Fix $x_0\in B_1$ and let $r$ be such that $B_r(x_0)\subset B_1$.
Let us denote by $a^+$ and $a^-$ the averages of $A$ in $H\cap B_r(x_0)$ and $H^c\cap B_r(x_0)$ respectively. In an analogous way we define $g^+$ and $g^-$ the averages of $G$ in $H\cap B_r(x_0)$ and $H^c\cap B_r(x_0)$. For $x\in B_r(x_0)$ we define
$$
\overline{A}:= a^+\mathbbm{1}_H+a^-\mathbbm{1}_{H^c},\qquad \qquad \overline{G}:= g^+\mathbbm{1}_H+g^-\mathbbm{1}_{H^c}.
$$
Notice that by assumption
\begin{equation}\label{MVhoelder}
|A(x)-\overline{A}(x)|\leq C_A r^{\sigma}\qquad\mbox{and}\qquad |G(x)-\overline{G}(x)|\leq C_G r^{\sigma}.
\end{equation}
Let $w$ be the solution of
\begin{equation*}
\begin{cases}
-\textnormal{div}(\overline{A}\D w)=\textnormal{div}\,\overline{G}\quad& \text{in }B_r(x_0),\\
w=v &\text{on }\partial B_r(x_0).
\end{cases}
\end{equation*}
The last equation can be rewritten as
\begin{equation}
\label{coeff cost}
\begin{cases}
-\textnormal{div}(a^+\D w^+)=0 & \text{in }B_r(x_0)\cap H,\\
-\textnormal{div}(a^-\D w^-)=0 & \text{in }B_r(x_0)\cap H^c,\\
w^+=w^-& \text{on }B_r(x_0)\cap\dd H,\\
\langle a^+\D w^+, e_n\rangle-\langle a^-\D w^-, e_n\rangle=\langle g^+, e_n\rangle-\langle g^-, e_n\rangle,&\text{on }B_r(x_0)\cap\dd H,\\
\end{cases}
\end{equation}
where $w^+:=w\mathbbm{1}_{B_r(x_0)\cap H}$, $w^-:=w\mathbbm{1}_{B_r(x_0)\cap H^c}$.
Set
\begin{equation*}
\overline{D}_c w:=\sum_{i=1}^n \overline{A}_{in}\D_i w+\langle\overline{G}, e_n\rangle,
\end{equation*}
where $\overline{A}_{in}$ is the $(i,n)$-th entry of the matrix $\overline{A}$. We notice that $\overline{D}_c w$ has no jumps on the boundary thanks to the transmission condition in \eqref{coeff cost}. 
This allows us to prove that the distributional gradient of $\overline{D}_c w$ coincides with the point-wise one.\\
\textbf{Step 1:} 
{\em Tangential derivatives of $w$}. Let us denote with $\tau$ the general direction tangent to the hyperplane $\partial H$. Since $\overline{A}$ and $\overline{G}$ are both constant along the tangential directions, the classical difference quotient method gives that $\D_{\tau}w\in W^{1,2}_{loc}(B_r(x_0))$ and 
$$\mbox{div}(\overline{A}\D (\D_{\tau}w))=0 \qquad\mbox{ in }  B_r(x_0).$$
Hence, Caccioppoli's inequality holds:
\begin{equation}
\label{Caccioppoli}
\int_{B_{\rho}(x)}|\D(\D_{\tau} w)|^2\,dy\leq \frac{c(n,\nu,N)}{\rho^2}\int_{B_{2\rho}(x)}|\D_{\tau} w-(\D_{\tau} w)_{x,2\rho}|^2\,dy,
\end{equation}
for all balls $B_{2\rho}(x)\subset B_r(x_0)$ and, by De Giorgi's regularity theorem, $\D_{\tau} w$ is H\"older continuous and there exists $\gamma=\gamma(n,\nu,N)>0$ such that if $B_{s}(x)\subset B_r(x_0)$ 
\begin{align}
\label{Eqn1}
\int_{B_{\rho}(x)}|\D_{\tau} w-(\D_{\tau} w)_{x,\rho}|^2\,dy \leq c(n,\nu,N)\bigg( \frac{\rho}{s} \bigg)^{n+2\gamma}\int_{B_{s}(x)}|\D_{\tau} w-(\D_{\tau} w)_{x,s}|^2\,dy,
\end{align}
for any $\rho\in \big(0,\frac{s}{2}\big)$ and
\begin{equation}
\label{Eqn2}
\max_{B_{\frac{\rho}{2}}(x)}|\D_{\tau} w|^2\leq \frac{c(n,\nu,N)}{\rho^n}\int_{B_{\rho}(x)}|\D_{\tau} w|^2\,dy.
\end{equation}
\textbf{Step 2:} {\em Regularity of $\overline{D}_c w$}. First of all observe that $\D_{\tau}(\overline{D}_c w)=\overline{D}_c(\D_{\tau} w)-\langle\overline{G}, e_n\rangle$. This implies by Step 1 that the tangential derivatives of $\overline{D}_c w$ belong to $L^2_{loc}(B_r(x_0))$. Furthermore we can estimate directly by definition of $\overline{D}_c w$:
\begin{equation*}
|\D_n(\overline{D}_c w)|\leq c(n,N)|\D \D_{\tau} w|,
\end{equation*}
which implies again by Step 1
\begin{equation*}
|\D\overline{D}_c w|\leq c(n,N)|\D \D_{\tau} w|.
\end{equation*}
We can conclude that $\overline{D}_c w\in W^{1,2}_{loc}(B_r(x_0))$. Using Poincaré's inequality and \eqref{Caccioppoli}, we have 
\begin{align*}
\int_{B_\rho(x)}|\overline{D}_c w-(\overline{D}_c w)_{x,\rho}|^2\,dy
& \leq c(n)\rho^2\int_{B_\rho(x)}|\D(\overline{D}_c w)|^2\,dy\\
& \leq c(n,N)\rho^2\int_{B_\rho(x)}|\D(\D_{\tau}w)|^2\,dy\\
& \leq c(n,\nu,N)\int_{B_{2\rho}(x)}|\D_{\tau}w-(\D_{\tau}w)_{x,2\rho}|^2\,dy,
\end{align*}
for any $B_{2\rho}(x)\subset B_r(x_0)$. By \eqref{Eqn1} we infer
\begin{equation*}
\begin{split}
&\int_{B_\rho(x)}|\overline{D}_c w-(\overline{D}_c w)_{x,\rho}|^2\,dy\\
& \leq c(n,\nu,N)\bigg( \frac{\rho}{r} \bigg)^{n+2\gamma}\int_{B_{\frac{r}{2}}(x)}|\D_{\tau} w-(\D_{\tau} w)_{x,\frac r2}|^2\,dy\\
&\leq c(n,\nu,N)\bigg( \frac{\rho}{r} \bigg)^{n+2\gamma}\int_{B_r(x_0)}|\D_{\tau} w|^2\,dy,
\end{split}
\end{equation*}
for any $x\in B_{\frac{r}{4}}(x_0)$ and $\rho\leq \frac{r}{4}$. Hence by Lemma 4.2 in \cite{MV} (see also \cite[Lemma 7.51]{AFP}), $\overline{D}_c w$ is H\"older continuous and by \eqref{Eqn2} we get:
\begin{equation}
\label{Eqn3}
\begin{split}
\max_{B_{\frac r4}(x_0)}|\overline{D}_c w|^2
& \leq c(n,\nu,N)\int_{B_r(x_0)}|\D_{\tau} w|^2\,dy+\bigg| \fint_{B_{\frac{r}{4}}(x_0)}\overline{D}_c w(y)\,dy \bigg|^2\\
& \leq \frac{c(n,\nu,N)}{r^n}\int_{B_r(x_0)}|\D w|^2\,dy+2\norm{G}^2_{L^\infty}.
\end{split}
\end{equation}
\textbf{Step 3:} {\em Comparison between $v$ and $w$.} Subtracting the equation for $w$ from the equation for $v$ we get
\begin{align}
& \int_{B_r(x_0)}\overline{A}_{ij}(x)\big(\D_i v-\D_i w\big)\D_j \varphi\,dx\notag\\
& =\int_{B_r(x_0)}\bigl(\overline{A}_{ij}(x)-A_{ij}(x)\bigr)\D_i v \D_j \varphi\,dx +\int_{B_r(x_0)}\bigl(\overline{G}_i-G_i\bigr)\D_i\varphi\,dx
\end{align}
for any $\varphi \in W^{1,2}_0(B_r(x_0))$. Choosing $\varphi =v-w$ in the previous equation 
and using assumption $\eqref{MVhoelder}$ we have
\begin{equation}
\label{Eqn6}
\nu\int_{B_r(x_0)}|\D v-\D w|^2\,dx\leq C_Ar^{\sigma}\int_{B_r(x_0)}|\D v|^2\,dy+C_G r^{n+\sigma}.
\end{equation}
Finally we can estimate
\begin{equation*}
\begin{split}
&\int_{B_\rho(x_0)}|\D v|^2\,dy\leq 2\int_{B_\rho(x_0)}|\D w|^2\,dy+2\int_{B_\rho(x_0)}|\D v-\D w|^2\,dy\\
& \leq 2\omega_n\rho^n \sup_{B_{\frac{r}{4}}}|\D w|^2+2\int_{B_\rho(x_0)}|\D v-\D w|^2\,dy,
\end{split}
\end{equation*}
for any $\rho\leq\frac{r}{4}$, and observing that 
\begin{equation*}
\begin{split}
\sup_{B_{\frac{r}{4}}(x_0)}|\D w|^2
&=\sup_{B_{\frac{r}{4}}(x_0)}|\D_{\tau} w|^2+\sup_{B_{\frac{r}{4}}(x_0)}|\D_n w|^2\\
& \leq c(n,\nu,N)\sup_{B_{\frac{r}{4}}(x_0)}|\D_{\tau} w|^2+c(\nu)\sup_{B_{\frac{r}{4}}(x_0)}|\overline{D}_c w|^2+c(\nu,\norm{G}_{\infty}),
\end{split}
\end{equation*}
by \eqref{Eqn2}, \eqref{Eqn3}, the minimality of $w$ and Young's inequality we gain
\begin{equation*}
\begin{split}
& \int_{B_\rho(x_0)}|\D v|^2\,dy\\
& \leq c(n,\nu,N)\bigg( \frac{\rho}{r}\bigg)^n\int_{B_r(x_0)}|\D w|^2\,dy+c(n,\nu,\norm{G}_{\infty},C_A,C_G)\bigg[r^{\sigma}\int_{B_r(x_0)}|\D v|^2\,dy+r^n\bigg]\\
& \leq C(n,\nu,N,\norm{G}_{\infty},C_A,C_G)\Bigg\lbrace\Bigg[ \bigg( \frac{\rho}{r}\bigg)^n+r^{\sigma} \Bigg]\int_{B_r(x_0)}|\D v|^2\,dy+r^n\Bigg\rbrace,
\end{split}
\end{equation*}
which leads to our aim if we apply Lemma \ref{Lemma iterativo 2}.
\end{proof}
The next lemma is inspired by \cite[Proposition 2.4]{FJ} and is the main result of this section.\\
In the sequel we shall consider the worst H\"older exponent introduced in  \eqref{Hoelderianity1} and \eqref{Hoelderianity2}, defined as
\begin{equation}\label{worst}
\delta:=\min{\left\{\alpha,\beta\right\}}.
\end{equation}
\begin{lemma}
\label{Lemma decadimento 1}
Let $(E,u)$ be a {$(\Lambda,\alpha)$}-minimizer of the functional ${\mathcal F}$ defined in \eqref{intro0}. There exists $\tau_0\in(0,1)$ such that the following statement is true: for all $\tau \in (0,\tau_0)$ there exists $\varepsilon_0=\varepsilon_0(\tau)>0$ such that if $B_r(x_0)\subset \subset \Omega$ with $r^{\frac \delta{2n}}<\tau$ and one of the following conditions holds:
\begin{itemize}
\item[\emph{(i)}]$ |E\cap B_r(x_0)|<\varepsilon_0 |B_r(x_0)|$,
\item[\emph{(ii)}]$ |B_r(x_0)\setminus E|<\varepsilon_0 |B_r(x_0)|$,
\item[\emph{(iii)}] There exists a halfspace $H$ such that $\frac{\left|(E\Delta H)\cap B_r(x_0)\right|}{|B_r(x_0)|}<\varepsilon_0$,
\end{itemize}
then
$$
\int_{B_{\tau r}(x_0)}|\nabla u|^2\,dx\leq C_4 \bigg[\tau ^{n} \int_{B_r(x_0)}|\nabla u|^2\,dx+r^{n}\bigg],
$$
for some positive constant $C_4=C_4\big(n,\nu,N,L,\alpha,\beta,L_\alpha,L_\beta,\norm{\D u}_{L^2(\Omega)}\big)$.
\end{lemma}
\begin{proof}
Let us fix $B_r(x_0)\subset \subset \Omega$ and $0<\tau<1$. Without loss of generality, we may assume that $\tau < 1/4$ and $x_0=0$. We start proving the assertion in the case (i), being the proof in the case (ii) similar.
Let us define
\begin{equation*}
A^0_{ij}:=a_{ij}(x_0,u_{r/2}(x_0)),\quad B^0_i:=a_i(x_0,u_{r/2}(x_0)),\quad  f^0:=a(x_0,u_{r/2}(x_0)),
\end{equation*}
\begin{equation*}
F_0(z):=\langle A^0 z,z\rangle+\langle B^0,z\rangle+f^0.
\end{equation*}
Let us denote by $v$ the solution of the following problem:
$$
\min_{w\in u+H^1_0(B_{r/2})}
\mathcal{F}_0(w;B_{r/2}),
$$
where
\begin{equation*}
\mathcal{F}_0(w;B_{r/2}):=\int_{B_{r/2}}F_0(\nabla w)\,dx.
\end{equation*}
Now we use the following identity 
\begin{equation*}
\langle A^0 \xi,\xi\rangle-\langle A^0 \eta,\eta\rangle=\langle A^0 (\xi-\eta),\xi-\eta\rangle+2\langle A^0 \eta,\xi-\eta\rangle, \quad\forall\xi,\eta\in\R^n,
\end{equation*}
in order to deduce that
\begin{align}\label{EF0}
& \mathcal{F}_0(u)-\mathcal{F}_0(v)\notag\\
& =\int_{B_{r/2}}\bigl[\langle A^0 \D u,\D u\rangle-\langle A^0 \D v,\D v\rangle\bigr]\,dx+\int_{B_{r/2}}\langle B^0,\D u-\D v\rangle\,dx\notag\\
& =\int_{B_{r/2}}\langle A^0 (\D u-\D v),\D u-\D v\rangle\,dx\notag\\
& +2\int_{B_{r/2}}\langle A^0 \D v,\D u-\D v\rangle\,dx + 
\int_{B_{r/2}}\langle B^0,\D u-\D v\rangle \,dx.
\end{align}
By the Euler-Lagrange equation for $v$ we deduce that the sum of the last two integrals in the previous identity is zero, being also $u=v$ on $\partial B_{r/2}$. Therefore, using the ellipticity assumption of $A^0$ we finally achieve that
\begin{equation}\label{EllF0}
\nu\int_{B_{r/2}}|\D u - \D v|^2\,dx\leq\mathcal{F}_0(u)-\mathcal{F}_0(v).
\end{equation}
Now we prove that $u$ is an $\omega$-minimizer of $\mathcal{F}_0$. We start writing
\begin{align}
\label{omegamin}
\mathcal{F}_0(u)
& =\mathcal{F}(E,u)+[\mathcal{F}_0(u)-\mathcal{F}(E,u)]\notag\\
&\leq \mathcal{F}(E,v)+[\mathcal{F}_0(u)-\mathcal{F}(E,u)]\notag\\
& = \mathcal{F}_0(v)+[\mathcal{F}_0(u)-\mathcal{F}(E,u)]+[\mathcal{F}(E,v)-\mathcal{F}_0(v)].
\end{align}
\emph{Estimate of $\mathcal{F}_0(u)-\mathcal{F}(E,u)$}. We use \eqref{Hoelderianity1}, \eqref{Hoelderianity2}, \eqref{ellipticity1}, \eqref{ellipticity2} and \eqref{hoelderu} to infer
\begin{align}
\label{eqq6}
& \mathcal{F}_0(u)-\mathcal{F}(E,u)=\int_{B_{r/2}}
\bigl(a_{ij}(x_0,u_{r/2}(x_0))-a_{ij}(x,u(x))\bigr)\D_iu\D_ju\,dx\notag\\
& +\int_{B_{r/2}}\bigl(a_i(x_0,u_{r/2}(x_0))-a_i(x,u(x))\bigr)\D_iu\,dx\notag\\
& +\int_{B_{r/2}}\bigl(a(x_0,u_{r/2}(x_0))-a(x,u(x))\bigr)\,dx -\int_{B_{r/2}\cap E}G(x,u,\D u)\,dx\notag\\
& \leq c\big(n,L_\alpha,L_{\beta}\norm{\D u}_{L^2(\Omega)}\big)\biggl(r^{\frac \delta 2}\int_{B_{r/2}}|\D u|^2\,dx+r^{n+\frac \delta 2}\biggr)+C(N,L)\bigg(\int_{B_{r/2}\cap E}|\D u|^2\,dx+r^n\bigg),
\end{align}
where we denoted $L_\alpha,L_\beta$ the greatest modulus of H\"older continuity of the data 
$a_{ij},b_{ij}, a_i,b_i,a,b$ defined in \eqref{Hoelderianity1} and \eqref{Hoelderianity2}.
Now we use H\"older's inequality and Lemma \ref{Lemma maggiore sommabilità} to estimate
\begin{align}\label{improvement}
\int_{B_{r/2}\cap E}|\D u|^2\,dx
& \leq |E\cap B_r|^{1-1/s}|B_r|^{1/s}\biggl(\fint_{B_{r/2}}|\D u|^{2s}\biggr)^{1/s}\notag\\
&\leq C_1^{1/s}\biggl(\frac{|E\cap B_r|}{|B_r|}\biggr)^{1-1/s}\int_{B_r}\big(1+|\D u|^2\big)\,dx.
\end{align}
Merging the last estimate in $\eqref{eqq6}$ we deduce
\begin{align}\label{E0}
\mathcal{F}_0(u)-\mathcal{F}(E,u)
& \leq\Bigl(c\big(n,L_\alpha,L_\beta,\norm{\D u}_{L^2(\Omega)}\big)
+C(N,L)C_1^{1/s}\Bigr)
\Bigl(r^{\frac \delta 2}+\varepsilon_0^{1-1/s}\Bigr)
\int_{B_r}|\D u|^2\,dx\nonumber\\
&+ \Bigl({C(N,L)C_{1}^{1/s}+C(N,L)}+c\big(n,L_\alpha,L_\beta,\norm{\D u}_{L^2(\Omega)}\big)\Bigr)r^n.
\end{align}
{\em Estimate of $\mathcal{F}(E,v)-\mathcal{F}_0(v)$}.
\begin{align}\label{E1}
\mathcal{F}(E,v)-\mathcal{F}_0(v)
& =\int_{B_{r/2}}
\bigl(a_{ij}(x,v(x))-a_{ij}(x_0,u_{r/2}(x_0))\bigr)\D_iv\D_jv\,dx\notag\\
&+\int_{B_{r/2}}\bigl(a_i(x,v(x))-a_i(x_0,u_{r/2}(x_0))\bigr)\D_iv\,dx \notag\\
& +\int_{B_{r/2}}\bigl(a(x,v(x))-a(x_0,u_{r/2}(x_0))\bigr)\,dx
+\int_{B_{r/2}\cap E}G(x,v,\D v)\,dx.
\end{align}
If we choose now $z\in\partial B_{r/2}$, recalling that $u(z)=v(z)$ we deduce
\begin{align*}
& \bigl|a_{ij}(x,v(x))-a_{ij}(x_0,u_{r/2}(x_0))\bigr|\\
& =\bigl|a_{ij}(x,v(x))-
a_{ij}(x,v(z))+a_{ij}(x,u(z))-a_{ij}(x_0,u_{r/2}(x_0))\bigr|\\
& \leq \big(L_\beta|v(x)-v(z)|^\beta+C\big(L_\beta,\norm{\nabla u}_{L^2(\Omega)}\big)r^{\frac \delta 2}+L_\alpha r^\delta\big)\\
& \leq \big(c(\beta,L_\beta)\text{osc}(u,\partial B_{r/2})^\beta+C(n,\nu,N,L,\beta,L_\alpha,L_\beta)r^\beta+C\big(L_\beta,\norm{\nabla u}_{L^2(\Omega)}\big)r^{\frac \delta 2}+r^\delta\big)\\
& \leq C\big(n,\nu,N,L,\beta,L_\alpha,L_\beta,\norm{\D u}_{L^2(\Omega)}\big)r^{\frac \delta 2},
\end{align*}
where we used the fact that $\text{osc}(v,B_{r/2})\leq \text{osc}(u,\partial B_{r/2})+C(n,\nu,N,L)r$  (see \cite[Lemma 8.4]{Giu}). Analogously we can estimate the other differences in $\eqref{E1}$, deducing
\begin{align*}
& \mathcal{F}(E,v)-\mathcal{F}_0(v)\leq C\big(n,\nu,N,L,\alpha,\beta,L_\alpha,L_\beta,\norm{\D u}_{L^2(\Omega)}\big) r^{\frac \delta 2}\biggl(\int_{B_{r/2}}|\D v|^2\,dx+ r^n\biggr)\\
& +C(N,L)\biggl(\int_{B_{r/2}\cap E}|\D v|^2\,dx +r^n\biggr),
\end{align*}
Reasoning in a similar way as in \eqref{improvement}, we can apply the higher integrability for $v$ given by Lemma \ref{Lemma maggiore sommabilità funzione confronto} and infer
$$
\int_{B_{r/2}\cap E}|\D v|^2\,dx \leq C(n,\nu,N,L)\varepsilon_0^{1-1/s}\biggl(\int_{B_{r}}|\D u|^2\,dx+r^n\biggr).
$$
Therefore we obtain
\begin{align}\label{E2}
\mathcal{F}(E,v)-\mathcal{F}_0(v) \leq C(n,\nu,N,L,\alpha,\beta,L_\alpha,L_\beta,\norm{\D u}_{L^2(\Omega)})\biggl[\Bigl(r^{\frac \delta 2}+\varepsilon_0^{1-1/s}\Bigr)\int_{B_r}|\D u|^2\,dx + r^n\biggr].
\end{align}
Finally, collecting $\eqref{EllF0}$, $\eqref{omegamin}$, $\eqref{E0}$ and $\eqref{E2}$, if we choose $\varepsilon_0$ such that $\varepsilon_0^{1-\frac{1}{s}}=\tau^n$, recalling that $r^{\frac{\delta}{2n}}<\tau$, we conclude that
\begin{equation}\label{FEstimate}
\int_{B_{r/2}}|\D u - \D v|^2\,dx\leq C \bigg[\tau^n\int_{B_r}|\D u|^2\,dx+ r^n\bigg],
\end{equation}
for some constant $C=C\big(n,\nu,N,L,\alpha,\beta,L_\alpha,L_\beta\norm{\D u}_{L^2(\Omega)}\big)$. On the other hand $v$ is the solution of a uniformly elliptic equation with constant coefficients, so we have
\begin{equation}\label{FEstimate1}
\int_{B_{\tau r}}|\D v|^2\,dx\leq C(n,\nu,N)\tau^n \int_{B_{r/2}}|\D v|^2\,dx\leq C(n,\nu,N,L)\bigg[\tau^n \int_{B_{r/2}}|\D u|^2\,dx+r^n\bigg].
\end{equation}
Hence we may estimate, using \eqref{FEstimate} and \eqref{FEstimate1},
\begin{align*}
\int_{B_{\tau r}}|\D u|^2\,dx\leq 2 \int_{B_{\tau r}}|\D v -\D u|^2\,dx +2 \int_{B_{\tau r}}|\D v|^2\,dx\leq C\bigg[\tau^n\int_{B_r}|\D u|^2\,dx+r^n\bigg],
\end{align*}
for some constant $C=C\big(n,\nu,N,L,\alpha,\beta,L_\alpha,L_\beta\norm{\D u}_{L^2(\Omega)}\big)$.\\
We are left with the case (iii). Let $H$ be the half-space from our assumption and let us denote accordingly
\begin{align*}
& A^0_{ij}(x):=a_{ij}(x,u(x))+\mathbbm{1}_{H}b_{ij}(x,u(x)) ,\\
& B^0_{ij}(x):=a_i(x,u(x))+\mathbbm{1}_{H}b_i(x,u(x)),\\
&f^0(x):= a(x,u(x))+\mathbbm{1}_{H} b(x,u(x)),\\
&F_0(x,z):=\langle A^0(x) z,z\rangle+\langle B^0(x),z\rangle+f^0(x).
\end{align*}
Let us denote by $v_H$ the solution of the following problem
$$
\min_{w\in u+H^1_0(B_{r/2})}
\mathcal{F}_0(w;B_{r/2}),
$$
where
\begin{equation*}
\mathcal{F}_0(w;B_{r/2}):=\int_{B_{r/2}}F_0(x,\nabla w)\,dx.
\end{equation*}
Let us point out that $v_H$ solves the Euler-Lagrange equation
\begin{equation}\label{MV}
-2\,\mathrm{div}(A^0\D v_H)=\textnormal{div}\,B^0 \quad \text{in }\mathcal{D}'(B_{r/2}).
\end{equation}
Therefore we are in position to apply Lemma \ref{Lemma semispazio} to the function $v_H$. Indeed, from the H\"older continuity of $u$ (see Remark \ref{Osservazione Holderianita}) we deduce that the restrictions of $A^0$ and $B^0$ onto $H\cap B_r$ and $B_r\setminus H$ respectively are H\"older continuous. We can conclude using also $\eqref{hoelderu}$ that there exist two constants $C=C\big(n,\nu,N,L,\alpha,\beta,L_\alpha,L_\beta\norm{\D u}_{L^2(\Omega)}\big)$ and $\tau_0=\tau_0
\big(n,\nu,N,L,\alpha,\beta,L_\alpha,L_\beta\norm{\D u}_{L^2(\Omega)}\big)$ such that for $\tau<\tau_0$
\begin{equation}\label{Decayvh}
\int_{B_{\tau r}}
|\D v_H|^2\,dx\leq 
C\bigg[\tau^n\int_{B_{r/2}}
|\D v_H|^2\,dx + r^n\bigg].
\end{equation}
In addition, using the ellipticity condition of $A^0$ we can argue as in $\eqref{EF0}$ to deduce using also the fact that $v_H$ satisfies $\eqref{MV}$,
\begin{equation}\label{EllF1}
\nu\int_{B_{r/2}}|\D u - \D v_H|^2\,dx\leq\mathcal{F}_0(u)-\mathcal{F}_0(v_H).
\end{equation}
One more time we can prove that $u$ is an $\omega$-minimizer of $\mathcal{F}_0$. We start as above writing
\begin{align}
\label{omegaminH}
&\mathcal{F}_0(u)
 =\mathcal{F}(E,u)+[\mathcal{F}_0(u)-\mathcal{F}(E,u)]\notag\\
&\leq \mathcal{F}(E,v_H)+[\mathcal{F}_0(u)-\mathcal{F}(E,u)]\notag\\
& = \mathcal{F}_0(v_H)+[\mathcal{F}_0(u)-\mathcal{F}(E,u)]+[\mathcal{F}(E,v_H)-\mathcal{F}_0(v_H)].
\end{align}
We can estimate the differences $\mathcal{F}_0(u)-\mathcal{F}(E,u)$ and $\mathcal{F}(E,v_H)-\mathcal{F}_0(v_H)$ exactly as before using this time the higher integrability given in Lemma \ref{Lemma maggiore sommabilità funzione confronto}.
We conclude that
\begin{equation}\label{FEHstimate}
\int_{B_{r/2}}|\D u - \D v_H|^2\,dx\leq C \bigg[\tau^n\int_{B_r}|\D u|^2\,dx+ r^n\bigg],
\end{equation}
for some constant $C=C\big(n,\nu,N,L,\alpha,\beta,L_\alpha,L_\beta\norm{\D u}_{L^2(\Omega)}\big)$. From the last estimate we can conclude the proof as before using \eqref{Decayvh} and \eqref{EllF1}.
\end{proof}
\bigskip
\section{Energy density estimates}
\label{Energy density estimates}
This section is devoted to prove a lower bound estimate for the functional $ {\mathcal F}(E,u;B_r (x_0))$. {Points i) and ii) of} Lemma \ref{Lemma decadimento 1} are the main tools to achieve such result. We shall prove that the energy $\mathcal F$ decays ``fast'' if the perimeter of $E$ is ``small''. In this section we will use a scaling argument.

\begin{lemma}[Scaling of $(\Lambda,\alpha)$-minimizers]
\label{Lemma riscalamento}
Let $B_r(x_0)\subset\Omega$ and let $(E,u)$ be a $(\Lambda,\alpha)$-minimizer of $\mathcal{F}$ in $B_r(x_0)$. Then $(E_{r},u_{r})$ is a $(\Lambda r^{\gamma},\alpha)$-minimizer of $\mathcal{F}_r$ in $B_1$, { for $\gamma=1+n(\alpha-1)\in(0,1)$} where
\begin{equation*}
E_{r}:=\frac{E-x_0}{r},\qquad u_{r}(y):=r^{-\frac{1}{2}}u(x_0+ry),\quad\text{for } y\in B_1,
\end{equation*}
\begin{align}
\label{RiscalamentoF}
\mathcal{F}_r(E_{r},u_{r};B_1):=
r\int_{B_1}\big[F(x_0+ry,r^{\frac{1}{2}}u_{r},r^{-\frac{1}{2}}\D u_{r})+\mathbbm{1}_{E_{r}}G(x_0+ry,r^{\frac{1}{2}}u_{r},r^{-\frac{1}{2}}\D u_{r}) \big]\,dy+P(E_{r};B_1).
\end{align}
\end{lemma}
\begin{proof}
Since $\D u_{r}(y)=r^{\frac{1}{2}}\D u(x_0+ry)$, for any $y\in B_1$, we rescale:
\begin{align*}
&\mathcal{F}(E,u;B_r(x_0))=r^n\int_{B_1}\big[F(x_0+ry,u(x_0+ry),\D u(x_0+ry))\\
&+\mathbbm{1}_{E}(x_0+ry)G(x_0+ry,u(x_0+ry),\D u(x_0+ry)) \big]\,dy+r^{n-1}P(E_{r};B_1)\\
& =r^{n-1}\mathcal{F}_r(E_{r},u_{r};B_1).
\end{align*}
Thus, if $\tilde{F}\subset\R^n$ is a set of finite perimeter with $\tilde{F}\Delta E_{r}\subset\subset B_1$ and $\tilde{v}\in H^1(B_1)$ is such that $\tilde{v}-u_{r}\in H_0^1(B_1)$, then
\begin{align*}
\mathcal{F}_r(E_{r},u_{r};B_1)
&=\frac{\mathcal{F}(E,u;B_r(x_0))}{r^{n-1}}\leq\frac{\mathcal{F}(F,v;B_r(x_0))+\Lambda|F\Delta E|^\alpha}{r^{n-1}}\\
&=\mathcal{F}_r(\tilde{F},\tilde{v};B_1)+\Lambda r^{\gamma}|\tilde{F}\Delta E_{r}|^\alpha,
\end{align*}
where $F:=x_0+r\tilde{F}$ and $v(x)= r^{\frac{1}{2}}\tilde{v}\big( \frac{x-x_0}{r} \big)$, for $x\in B_r(x_0)$.
\end{proof}

\begin{lemma}
\label{Lemma decadimento 2} Let $(E,u)$ be a $(\Lambda,\alpha)$-minimizer in $\Omega$ of the functional ${\mathcal F}$ defined in \eqref{intro0}. For every $\tau\in (0,1)$ there exists $\varepsilon_1=\varepsilon_1(\tau)>0$ such that, if $B_r(x_0)\subset \Omega$ and $P(E;B_r(x_0))<\varepsilon_1 r^{n-1}$, then
\begin{equation}\label{D}
\mathcal F(E,u;B_{\tau r}(x_0))\leq C_5 \bigl(\tau^n\mathcal F(E,u;B_r(x_0))+(\tau r)^{n\alpha}\bigr),
\end{equation}
for some positive constant $C_5=C_5\big(n,\nu,N,L,L_\alpha,L_\beta,\alpha,\beta,\Lambda,\norm{\D u}_{L^2(\Omega)}\big)$ independent of $\tau$ and $r$.
\end{lemma}
\begin{proof}
Let $\tau\in(0,1)$ and $B_r(x_0)\subset\Omega$. Without loss of generality, we may assume that $\tau<\frac 12$. We may also assume that $x_0=0$, and $r=1$ by scaling $E_r=\frac{E-x_0}{r}$, $u_r(y)=r^{-\frac 12}{u(x_0+ry)}$ for $y\in B_1$, and replacing $\Lambda$ with $\Lambda r^\gamma$. Thus, we have that $(E_r,u_r)$ is a $(\Lambda r^\gamma,\alpha)$-minimizer of $\mathcal F_r$ in $\frac{\Omega-x_0}{r}$. For simplicity of notation we can still denote $E_r$ by $E$, $u_r$ by $u$ and then, recalling that $\mathcal{F}=r^{n-1}\mathcal{F}_r$ and
$\gamma=n\alpha-(n-1)$, we have to prove that there exists $\varepsilon_1=\varepsilon_1(\tau)$ such that, if $P(E;B_1)<\varepsilon_1$, then
$$
\mathcal F_r(E,u;B_{\tau })\leq C_5 \big(\tau^n\mathcal F_r(E,u;B_1)+\tau^{n\alpha}r^\gamma\big).
$$
Note that, since $P(E;B_1)<\varepsilon_1$, by the relative isoperimetric inequality, either $|B_1\cap E|$ or $|B_1\setminus E|$ is small and thus Lemma \ref{Lemma decadimento 1} can be applied. Assuming that $|B_1\setminus E|\leq |B_1\cap E|$ and using the relative isoperimetric inequality we can deduce that
$$
|B_1\setminus E|\leq c(n)P(E;B_1)^{\frac{n}{n-1}}.
$$
If we choose as a representative of $E$ the set of points of density one, we get, by Fubini's theorem that
$$
|B_1\setminus E|\geq \int_{\tau}^{2\tau}\mathcal H^{n-1}(\partial B_\rho\setminus E)\,d\rho.
$$
Combining these inequalities, we can choose $\rho\in (\tau, 2\tau)$ such that
\begin{equation}\label{perimeter}
\mathcal{H}^{n-1}(\partial B_\rho\setminus E)\leq \frac{c(n)}{\tau}P(E;B_1)^{\frac{n}{n-1}}\leq \frac{c(n)\varepsilon_1^{\frac{1}{n-1}}}{\tau}P(E;B_1).
\end{equation}
Now we set $F=E\cup B_\rho$ and observe that
$$
P(F;B_1)\leq P(E;B_1\setminus \overline{B}_{\rho})+\mathcal H^{n-1}(\partial B_{\rho}\setminus E).
$$
If we choose $(F,u)$ to test the $(\Lambda r^\gamma,\alpha)$-minimality of $(E,u)$ we get
\begin{align*}
\mathcal{F}_r(E,u)
& \leq\mathcal{F}_r(F,u)+ \Lambda r^\gamma|F\setminus E|^\alpha\\
& \leq P(E;B_1\setminus \overline{B}_{\rho})+\mathcal H^{n-1}(\partial B_{\rho}\setminus E)+ \Lambda r^\gamma |B_{\rho}|^\alpha\\
& +r \int_{B_1}\big ( F(x_0+ry,r^{\frac 12} u(y),r^{-\frac{1}{2}}\nabla u(y))+\mathbbm{1}_{F}G(x_0+ry,r^{\frac 12} u(y),r^{-\frac{1}{2}}\nabla u(y))\big )\,dy .
\end{align*}
Then getting rid of the common terms we obtain
\begin{equation*}
P(E;B_\rho)\leq \mathcal H^{n-1}(\partial B_{\rho}\setminus E)+ r\int_{B_{\rho}} G(x_0+ry,r^{\frac 12} u(y),r^{-\frac{1}{2}}\nabla u(y))\,dy + \Lambda r^\gamma |B_{\rho}|^\alpha.
\end{equation*}
Now if we choose $\varepsilon_1$ such that $c(n)\varepsilon_1^{\frac{1}{n-1}}\leq \tau^{n+1}$ we have from \eqref{perimeter}
\begin{equation*}
P(E;B_\rho)\leq \tau^n P(E;B_1)+ r\int_{B_{\rho}} G(x_0+ry,r^{\frac 12} u(y),r^{-\frac{1}{2}}\nabla u(y))\,dy + \Lambda r^\gamma |B_{\rho}|^\alpha.
\end{equation*}
Then, we choose $\varepsilon_1$ satisfying $c(n)\varepsilon_1^{\frac{n}{n-1}} \leq \varepsilon_0(2\tau)|B_1|$ to obtain, using Lemma \ref{Lemma decadimento 1} and growth conditions \eqref{ellipticity1}, \eqref{ellipticity2},
\begin{align*}
r\int_{B_{\rho}} G(x_0+ry,r^{\frac 12} u(y),r^{-\frac{1}{2}}\nabla u(y))\,dy 
& \leq C(N,L)\int_{B_{\rho}} (|\nabla u|^2+r)\,dy\\
&\leq C\big(n,\nu,N,L,L_\alpha,L_{\beta},\alpha,\beta,\norm{\D u}_{L^2(\Omega)}\big) \tau^n\int_{B_{1}}
(|\nabla u|^2+r)\,dy.
\end{align*}
Finally, we recall that $\rho\in (\tau,2\tau)$ to conclude, using the previous estimates,

\begin{align*}
P(E;B_\tau)
& \leq C\big(n,\nu,N,L,L_\alpha,L_{\beta},\alpha,\beta,\norm{\D u}_{L^2(\Omega)}\big) \tau^n \bigg[\int_{B_{1}} (|\D u|^2+r)\,dy+P(E;B_1)\bigg]  + \Lambda r^\gamma |B_{2\tau}|^\alpha\\
& \leq C\big(n,\nu,N,L,L_\alpha,L_{\beta},\alpha,\beta,\norm{\D u}_{L^2(\Omega)}\big) \big[\tau^n\mathcal{F}_r(E,u;B_1)+\tau^{n\alpha}r^\gamma\big].
\end{align*}
From this estimate the result easily follows applying again Lemma \ref{Lemma decadimento 1}.
\end{proof}
In the sequel we will assume that the representative of the set $E$ is choosen in such a way that the topological boundary $\partial E$ concides with the closure of the reduced boundary, that is $\dd E=\overline{\dd\displaystyle^*E}$, (see also \cite{Ma} Proposition 12.19).
\begin{theorem}[Density lower bound]
\label{Density lower bound}
Let $(E,u)$ be a $(\Lambda,\alpha)$-minimizer of $\mathcal{F}$ in $\Omega$ and $U\subset \subset \Omega$ be an open set. Then there exists a constant $C_6=C_6\big(n,\nu,N,L,\alpha,\beta,L_\alpha,L_\beta,\Lambda,\norm{\D u}_{L^2(\Omega)},U\big)$,
such that, for every $x_0\in \partial E$ and $B_r(x_0)\subset U$, it holds
$$
P(E;B_r(x_0))\geq C_6r^{n-1}.
$$
Moreover, $\mathcal{H}^{n-1}((\partial E\setminus \partial^*E)\cap \Omega)=0$.
\end{theorem}

\begin{proof}
We start assuming that $x_0\in\displaystyle\dd^*E$. Without loss of generality we may also assume that $x_0=0$. Let 
\begin{equation*}
\tau\in\bigl(0,2^{-\frac{1}{\gamma}}\bigr)\,\text{ such that }\, 2C_5\tau^{n(1-\alpha)}<1,
\end{equation*}
\begin{equation*}
\sigma\in\bigl(0,1\bigr)\text{ such that }\, 2C_5C_3\sigma^\gamma<\varepsilon_1(\tau),\;\;2\omega_n\frac{L^2}{\nu}\sigma<\varepsilon_1(\tau),\;\;\sigma^\gamma<\tau^{n(1-\alpha)},
\end{equation*}
\begin{equation*}
0<r_0<\min\big\{1,C_3^{\frac{1}{\gamma}},\varepsilon_1(\tau)^{\frac{1}{\gamma}}\big\},
\end{equation*}
where $C_5$ and $\varepsilon_1$ come from Lemma \ref{Lemma decadimento 2}, $C_3$ comes from Theorem \ref{Energy upper bound}. 
We point out that $\tau,\sigma,r_0,\varepsilon_1(\sigma)$ depend on
$n,\nu,N,L,\alpha,\beta,L_\alpha,L_\beta,\Lambda,\norm{\D u}_{L^2(\Omega)}$ through the constants $C_3$ and $C_5$ only.
Let us suppose by contradiction that there exists $B_r\subset U$, with $r<r_0$, such that $P(E;B_r)<\varepsilon_1(\sigma)r^{n-1}$. We shall prove that
\begin{equation}
\label{Relazione iterativa}
\mathcal{F}(E,u;B_{\sigma\tau^hr})\leq\varepsilon_1(\tau)\tau^{\gamma h}(\sigma\tau^h r)^{n-1},
\end{equation}
for any $h\in\N_0$, reaching a contradiction afterward.\\
For $h=0$, using Lemma \ref{Lemma decadimento 2} with $\varepsilon_1=\varepsilon_1(\sigma)$, Theorem \ref{Energy upper bound}, $r<r_0< C_3^{\frac{1}{\gamma}}$ and $2C_5C_3\sigma^\gamma<\varepsilon_1(\tau)$, we get:
\begin{align*}
\mathcal{F}(E,u;B_{\sigma r})
&\leq C_5\big (\sigma^n\mathcal{F}(E,u;B_r)+(\sigma r)^{n\alpha}\big)\\
& \leq C_5C_3\sigma^nr^{n-1}+C_5\sigma^{n\alpha}r^{n-1}r^{\gamma}\\
& \leq 2C_5C_3\sigma^{n\alpha} r^{n-1}\leq\varepsilon_1(\tau)(\sigma r)^{n-1}.
\end{align*}
In order to prove the induction step we have to ensure to be in position to apply 
Lemma \ref{Lemma decadimento 2}, that is by proving smallness of the perimeter. In such regard, let us observe that, by the definition of $\mathcal{F}(E,u;B_\rho)$ and the growth condition given in \eqref{ggrowth},
$$
P(E;B_\rho)\leq\mathcal{F}(E,u;B_\rho)+2\omega_n\frac{L^2}{\nu}\rho^n,
$$  
for any $B_\rho\subset \Omega$.\\
Assuming that the induction hypothesis \eqref{Relazione iterativa} holds true for some $h\in\N$ and, being $2\omega_n\frac{L^2}{\nu}\sigma<{\varepsilon_1(\tau)}$,
$\tau<{2^{-\frac{1}{\gamma}}}$ and $r<1$, we infer
\begin{align*}
P(E;B_{\sigma\tau^{h}r})
& \leq\mathcal{F}(E,u;B_{\sigma\tau^{h}r})+2\omega_n\frac{L^2}{\nu}(\sigma\tau^{h}r)^n\\
&\leq (\sigma\tau^hr)^{n-1}\bigg( \varepsilon_1(\tau)\tau^{\gamma h}+2\omega_n
\frac{L^2}{\nu}\sigma\tau^h r\bigg)
\leq (\sigma\tau^hr)^{n-1}\varepsilon_1(\tau)(\tau^{\gamma h}+\tau^h)\\
& \leq (\sigma\tau^hr)^{n-1}\varepsilon_1(\tau)2\tau^\gamma\leq (\sigma\tau^hr)^{n-1}\varepsilon_1(\tau).
\end{align*}
We are now in position to apply Lemma \ref{Lemma decadimento 2} with $\varepsilon_1=\varepsilon_1(\tau)$. Using also the induction hypothesis and, since  $\sigma^\gamma<\tau^{n(1-\alpha)}$, $r<r_0\leq\varepsilon_1(\tau)^{\frac{1}{\gamma}}$  and $2C_5\tau^{n(1-\alpha)}<1$, we estimate:
\begin{align*}
\mathcal{F}(E,u;B_{\sigma\tau^{h+1}r})
& \leq C_5\big[ \tau^n\mathcal{F}(E,u;B_{\sigma\tau^hr})+\tau^{n\alpha}(\sigma\tau^hr)^{n\alpha} \big]\\
& \leq C_5\big[ \tau^n\varepsilon_1(\tau)\tau^{\gamma h}(\sigma\tau^hr)^{n-1}+\tau^{n\alpha}(\sigma\tau^hr)^{n\alpha} \big]\\
& =\tau^{\gamma h}(\sigma\tau^hr)^{n-1} C_5\big[ \tau^n\varepsilon_1(\tau)+\tau^{n\alpha}(\sigma r)^\gamma \big]\leq \tau^{\gamma h}(\sigma\tau^hr)^{n-1} \tau^{n}\big[ C_5\varepsilon_1(\tau)+C_5 r^\gamma \big]\\
& \leq \tau^{\gamma h}(\sigma\tau^hr)^{n-1}\tau^n 2C_5\varepsilon_1(\tau)\leq \tau^{\gamma h}(\sigma\tau^hr)^{n-1}\tau^n\varepsilon_1(\tau)\tau^{n(\alpha-1)}\\
& =\tau^{\gamma(h+1)}(\sigma\tau^{h+1}r)^{n-1}\varepsilon_1(\tau).
\end{align*}
We conclude that \eqref{Relazione iterativa} holds for any $h\in\N_0$. Thus, we gain
\begin{align*}
P(E;B_{\sigma\tau^hr})
& \leq\varepsilon_1(\tau)\tau^{\gamma h}(\sigma\tau^h r)^{n-1}+2\omega_n\frac{L^2}{\nu}(\sigma\tau^h r)^n\\
& \leq (\sigma\tau^hr)^{n-1}\tau^{\gamma h}\bigg( \varepsilon_1(\tau)+{2\omega_n\frac{L^2}{\nu}}\sigma\tau^{h(1-\gamma)} \bigg)\\
& \leq (\sigma\tau^hr)^{n-1}\tau^{\gamma h} \varepsilon_1(\tau)\big(1+\tau^{h(1-\gamma)} \big)\\
& \leq 2(\sigma\tau^hr)^{n-1}\tau^{\gamma h} \varepsilon_1(\tau).
\end{align*}
We finally get
\begin{equation*}
\lim_{\rho\rightarrow 0^+}\frac{P(E;B_\rho)}{\rho^{n-1}}=\lim_{h\rightarrow +\infty}\frac{P(E; B_{\sigma\tau^hr})}{(\sigma\tau^hr)^{n-1}}\leq\lim_{h\rightarrow+\infty}2\varepsilon_1(\tau)\tau^{\gamma h}=0,
\end{equation*}
which implies that $x_0\not\in\dd^*E$, that is a contradiction.
We recall that we chose the representative of $\dd E$ such that $\dd E=\overline{\dd\displaystyle^*E}$.
Thus, if $x_0\in\dd E$, there exists $(x_h)_{h\in\N}\subset\dd^*E$ such that $x_h\rightarrow x_0$ as $h\rightarrow+\infty$,
\begin{equation*}
P(E;B_r(x_h))\geq c\big(n,\nu,N,L,\alpha,\beta,L_\alpha,L_\beta,\Lambda,\norm{\D u}_{L^2(\Omega)}\big)r^{n-1}
\end{equation*}
and $B_r(x_h)\subset U$, for $h$ large enough. Passing to the limit as $h\rightarrow+\infty$, we get the thesis.
\end{proof}

\section{Compactness for sequences of minimizers}
\label{Compactness for sequences of minimizers}
In this section we basically follow the path given in \cite[Part III]{Ma}.
We start proving a standard compactness result.
\begin{lemma}[Compactness]
\label{Lemma compattezza}
Let $(E_h,u_h)$ be a sequence of $(\Lambda_h,\alpha)$-minimizers of $\mathcal{F}$ in $\Omega$ such that $\sup_h \mathcal F(E_h,u_h;\Omega)<\infty$ and 
$\Lambda_h\rightarrow \Lambda\in \mathbb R^+$. There exist a (not relabelled) subsequence and a $(\Lambda,\alpha)$-minimizer of $\mathcal F$ in $\Omega$,  $(E,u)$,  such that for every open set $U\subset \subset \Omega$, it holds
$$
E_h\rightarrow E \mbox { in } L^1(U),\quad u_h\rightarrow u \mbox { in } H^{1}(U),\quad P(E_h;U)\rightarrow P(E;U).
$$
In addition,
\begin{align}
& \label{boundary1} \mbox{if }x_h\in \partial E_h\cap U \mbox{ and } x_h\rightarrow x \in U, \mbox { then } x\in \partial E \cap U,\\
& \label{boundary2}\mbox{if }x\in \partial E\cap U, \mbox{ there exists } x_h\in \partial E_h\cap U \mbox{ such that } x_h\rightarrow x.
\end{align}
Finally, if we assume also that $\nabla u_h \rightharpoonup 0$ weakly in $L^2_{loc}(\Omega,\mathbb R^n)$ and $\Lambda_h \rightarrow 0$, as $h\rightarrow \infty$, then $E$ is a local minimizer 
of the perimeter, that is
$$
P(E;B_r(x_0))\leq P(F;B_r(x_0)),
$$
for every set $F$ such that $F\Delta E \subset\subset B_r(x_0)\subset \Omega.$
\end{lemma}
\begin{proof}
We start observing that, by the boundedness condition on $\mathcal F(E_h,u_h;\Omega)$, we may assume that $u_h$ weakly converges to $u$ in $H^{1}(U)$ and strongly in $L^2(U)$, and 
$\mathbbm{1}_{E_h}$ converges to $\mathbbm{1}_{E}$ in $L^1(U)$, as $h\rightarrow \infty$. By lower semicontinuity we are going to prove the $(\Lambda,\alpha)$-minimality of $(E,u)$. Let us fix 
$B_r(x_0)\subset\subset\Omega$ and assume for simplicity of notation that $x_0=0$. Let $(F,v)$ be a test pair such that $F\Delta E \subset \subset B_r$ 
and supp$(u-v)\subset \subset B_r$. 
We can handle the perimeter term as in \cite{Ma}, that is, eventually passing to a subsequence and using Fubini's theorem, we may choose $\rho<r$ such that, once again, 
$F\Delta E \subset \subset B_{\rho}$ and supp$(u-v)\subset \subset B_{\rho}$, and, in addition,
$$
\mathcal{H}^{n-1}(\partial^*F\cap \partial B_{\rho})=\mathcal{H}^{n-1}(\partial^*E_h\cap \partial B_{\rho})=0,
$$
and
\begin{equation}
\label{mismatch}
\lim_{h\rightarrow 0}\mathcal{H}^{n-1}(\partial B_{\rho}\cap(E\Delta E_h))=0.
\end{equation}
Now we choose a cut-off function $\psi\in C_0^1(B_r)$ such that $\psi\equiv 1$ in $B_{\rho}$ and define $v_h=\psi v+(1-\psi)u_h$, 
$F_h:=(F\cap B_{\rho})\cup (E_h\setminus B_{\rho})$ to test the minimality of $(E_h,u_h)$. Thanks to the $(\Lambda_h,\alpha)$-minimality of $(E_h,u_h)$ we have
\begin{align}\label{min}
& \int_{B_r} \bigl ( F(x,u_h,\nabla u_h)+\mathbbm{1}_{E_h}G(x,u_h,\nabla u_h)\bigr )\,dx+ P(E_h;B_r)\leq\notag\\
& \leq \int_{B_r} \bigl ( F(x,v_h,\nabla v_h)+\mathbbm{1}_{F_h}G(x,v_h,\nabla v_h)\bigr )\,dx+P(F_h;B_r)
+\Lambda_h|F_h\Delta E_h|^\alpha \notag\\
& \leq \int_{B_r} \bigl (F(x,v_h,\nabla v_h)+\mathbbm{1}_{F_h}G(x,v_h,\nabla v_h)\bigr )\,dx +P(F;B_{\rho})+\Lambda_h|F_h\Delta E_h|^\alpha \notag\\
& +P(E_h;B_r\setminus \overline{B}_{\rho})+\varepsilon_h.
\end{align}
The mismatch term $\varepsilon_h=\mathcal{H}^{n-1}(\partial B_{\rho}\cap (F^{(1)}\Delta E_{h}^{(1)}))$ appears because $F$ is not in general a compact variation of $E_h$. Nevertheless we have that  $\varepsilon_h\rightarrow 0$ because of the assumption $\eqref{mismatch}$ (see also \cite[Theorem 21.14]{Ma}).\\
Now we use the convexity of $F$ and $G$ with respect to the $z$ variable to deduce
\begin{align*}
& \int_{B_r} \bigl ( F(x,v_h,\nabla v_h)+\mathbbm{1}_{F_h}G(x,v_h,\nabla v_h)\bigr )\,dx\\
& \leq \int_{B_r} \bigl ( F(x,v_h,\psi \nabla v+(1-\psi)\nabla u_h)+\mathbbm{1}_{F_h}G(x,v_h,\psi \nabla v+(1-\psi)\nabla u_h)\bigr )\,dx\\
& + \int_{B_r}  \left < \nabla_z F(x,v_h,\nabla v_h),\nabla \psi(v-u_h)\right >\,dx + \int_{B_r}  \mathbbm{1}_{F_h} \left < \nabla_z G(x,v_h,\nabla v_h),\nabla \psi(v-u_h)\right >\,dx,\\
\end{align*}
where the last two terms in the previous estimate tend to zero as $h\rightarrow \infty$. Indeed, the term $\nabla \psi(v-u_h)$ strongly converges to zero in $L^2$, being $u=v$ in $B_r\setminus B_{\rho}$ and the first part in the scalar product weakly converges in $L^2$.
Then using again the convexity of $F$ and $G$ with respect to the $z$ variable we obtain, for some infinitesimal $\sigma_h$,
\begin{align}\label{compact}
&\int_{B_r} \big ( F(x,v_h,\nabla v_h)+\mathbbm{1}_{F_h}G(x,v_h,\nabla v_h)\big)\,dx \notag\\
& \leq \int_{B_r} \psi \big( F(x,v_h,\nabla v  )+\mathbbm{1}_{F_h}G(x,v_h,\nabla v)\big)\,dx\notag\\ 
& + \int_{B_r} (1-\psi) \big( F(x,v_h,\nabla u_h )+\mathbbm{1}_{F_h}G(x,v_h,\nabla u_h)\big)\,dx+ \sigma_h.
\end{align}
Finally, we combine $\eqref{min}$ and $\eqref{compact}$ and pass to the limit as $h\rightarrow+\infty$, using the lower semicontinuity on the left-hand side. For the right-hand side we observe that $\mathbbm{1}_{E_h}\rightarrow \mathbbm{1}_{E}$  and $\mathbbm{1}_{F_h}\rightarrow \mathbbm{1}_{F}$ in $L^1(B_r)$ and we use also the equi-integrability of $\left\{\nabla u_h\right\}_h$ to conclude,
\begin{align*}
& \int_{B_r} \psi\big ( F(x,u,\nabla u)+\mathbbm{1}_{E}G(x,u,\nabla u)\big )\,dx+ P(E;B_{\rho})\\
& \leq \int_{B_r} \psi \big ( F(x,v,\nabla v)+\mathbbm{1}_{F}G(x,v,\nabla v)\big )\,dx+ P(F;B_{\rho})+ \Lambda|F\Delta E|^\alpha.
\end{align*}
Letting $\psi \downarrow \mathbbm{1}_{B_{\rho}}$ we finally get
\begin{align}\label{lmini}
& \int_{B_{\rho}} \big ( F(x,u,\nabla u)+\mathbbm{1}_{E}G(x,u,\nabla u)\big )\,dx+ P(E;B_{\rho})\notag\\
& \leq \int_{B_{\rho}}  \big ( F(x,v,\nabla v)+\mathbbm{1}_{F}G(x,v,\nabla v)\big )\,dx+ P(F;B_{\rho})+ \Lambda|F\Delta E|^\alpha,
\end{align}
and this proves the $(\Lambda,\alpha)$-minimality of $(E,u)$.
\\To prove the strong convergence of $\nabla u_{h}$ to $\nabla{u}$ in $L^2(B_r)$ we start observing that 
by $\eqref{min}$ and $\eqref{compact}$ applied using $(E_h,u)$ to test the $(\Lambda,\alpha)$-minimality of $(E_h,u_h)$ we get
\begin{align*}
\int_{B_r} \psi\big ( F(x,u_h,\nabla u_h)+\mathbbm{1}_{E_h}G(x,u_h,\nabla u_h)\big )\,dx\leq
\int_{B_r} \psi \big( F(x,u,\nabla u  )+\mathbbm{1}_{E_h}G(x,u,\nabla u)\big)\,dx+ \sigma_h.
\end{align*}
Then from the equi-integrability of $\left\{\nabla u_h\right\}_h$  in $L^2(U)$ and recalling that $\mathbbm{1}_{E_h}\rightarrow \mathbbm{1}_{E}$ in $L^1(U)$, we obtain
\begin{align*}
\limsup_{h\rightarrow \infty}\int_{B_r} \psi \big( F(x,u_h,\nabla u_h  )+\mathbbm{1}_{E_h}G(x,u_h,\nabla u_h)\big)\,dx\leq \int_{B_r} \psi \big( F(x,u,\nabla u  )+\mathbbm{1}_{E}G(x,u,\nabla u)\big)\,dx.
\end{align*}
The opposite inequality can be obtained by semicontinuity. Thus we get
\begin{align}\label{norm}
\lim_{h\rightarrow \infty}\int_{B_r} \psi \bigl( F(x,u_h,\nabla u_h  )+\mathbbm{1}_{E_h}G(x,u_h,\nabla u_h)\bigr)\,dx= \int_{B_r} \psi \bigl( F(x,u,\nabla u  )+\mathbbm{1}_{E}G(x,u,\nabla u)\bigr)\,dx.
\end{align}
From the ellipticity condition in $\eqref{ellipticity1}$ we infer, for some $\sigma_h\rightarrow 0$,
\begin{align}
\nu\int_{B_r}\psi|\nabla u_h-\nabla u|^2\,dx
& \leq \int_{B_r}\psi\bigl(F(x,u_h,\nabla u_h)-F(x,u,\nabla u)\bigr)\,dx \nonumber \\
&+\int_{B_r}\psi\mathbbm{1}_{E}\bigl(G(x,u_h,\nabla u_h)-G(x,u,\nabla u)\bigr)\,dx  + \sigma_h.
\end{align}
Passing to the limit we obtain
$$
\lim_{h\rightarrow \infty}\int_{B_r}\psi|\nabla u_h -\nabla u |^2\,dx=0.
$$
Finally testing the minimality of $(E_h,u_h)$ with respect to the pair $(E,u)$ we also get
$$
\lim_{h\rightarrow \infty} P(E_h;B_{\rho})=P(E;B_{\rho}).
$$
With a usual argument we can deduce $u_h\rightarrow u$ in $W^{1,2}(U)$ and $P(E_h;U)\rightarrow P(E;U)$, for every open set $U\subset \subset \Omega$.
The topological information stated in $\eqref{boundary1}$ and $\eqref{boundary2}$ follows as in \cite[Theorem 21.14]{Ma} because it does not depend on the presence of the integral bulk part.
\end{proof}

\section{Height bound and Lipschitz approximation}
\label{Height bound and Lipschitz approximation}
In the following for $R>0$ and  $\nu\in\mathbbm{S}^{n-1}$ we will denote
\begin{equation*}
\mathbf{C}_{R}(x_0,\nu):=x_0+\{y\in\R^n\,:\,|\langle y,\nu \rangle|<R,\,|y-\langle y,\nu\rangle \nu|<R\},
\end{equation*}
the cylinder centered in $x_0$ with radius $R$ oriented in the direction $\nu$.\\
The cylinder of radius $R$ oriented in the direction $e_n$ with height 2 will be denoted as
\begin{equation*}
\mathbf{K}_{R}(x_0):=\{y=(y',y_n)\in\R^n\,:\,|y'-x_0'|<R,\,|y_n-(x_0)_n|<1\},
\end{equation*}
In addition we introduce some usual quantities involved in regularity theory  
\begin{definition}
Let $E$ be a set of locally finite perimeter, $x\in\dd E$, $r>0$ and $\nu\in \mathbbm S^{n-1}$. We define:
\begin{itemize}
\item the \textbf{cylindrical excess} of $E$ at the point $x$, at the scale $r$ and with respect to the direction $\nu$, as
\begin{equation*}
\ecc(x,r,\nu):=\frac{1}{r^{n-1}}\int_{\textbf{C}(x,r,\nu)\cap\dd^*E}\frac{|\nu_E-\nu|^2}{2}d\,\mathcal{H}^{n-1}=\frac{1}{r^{n-1}}\int_{\textbf{C}(x,r,\nu)\cap\dd^*E}[1-\langle\nu_E,\nu\rangle]\,d\mathcal{H}^{n-1}.
\end{equation*}
\item the \textbf{spherical excess} of $E$ at the point $x$, at the scale $r$ and with respect to the direction $\nu$, as
\begin{equation*}
{\mathbf e}(x,r,\nu):= \frac{1}{r^{n-1}}\int_{\partial E\cap B_r(x)}\frac{|\nu_E-\nu|^2}{2}d\mathcal H^{n-1}.
\end{equation*}
\item the \textbf{spherical excess} of $E$ at the point $x$ and at the scale $r$, as
\begin{equation*}
{\mathbf e}(x,r):=\min_{\nu \in \mathbbm S^{n-1}}{\mathbf e}(x,r,\nu).
\end{equation*}
\end{itemize} 
\end{definition}
{In the following, for simplicity of notation we will denote 
$$
\mathbf{C}_R=\mathbf{C}_R(0,e_n)=\{y=(y',y_n)\in\R^n\,:\,|y'|<R,\,|y_n|<R\}.
$$}
The following height bound lemma is a standard step in the proof of regularity because it is one of the main ingredients to prove the Lipschitz approximation theorem. {The results contained in this section are a consequence of the compactness lemma, the density lower bound and the lower semicontinuity of the excess. In the statement of these results we assume that $(E,u)$ is a $(\Lambda,\alpha)$-minimizer of $\mathcal{F}$. However the minimality is not used except to ensure compactness and the density lower bound}. 
\begin{lemma}[Height bound]
\label{height bound}
Let $(E,u)$ be a $(\Lambda,\alpha)$-minimizer of $\mathcal{F}$ in $B_r(x_0)$. There exist two positive constants $\varepsilon_2$ and $C_7$, depending on $n,\nu,N,L,\alpha,\beta,L_\alpha,L_\beta,\Lambda,\norm{\D u}_{L^2(\Omega)}$, such that if $x_0\in \partial E$ and
$$
{\mathbf e}(x,r,\nu)<\varepsilon_2,
$$
for some $\nu\in \mathbbm S^{n-1}$, then
$$
\sup_{y\in \partial E\cap B_{r/2}(x_0)}\frac{|\langle \nu,y-x_0\rangle|}{r}\leq C_7{\mathbf e}(x,r,\nu)^{\frac{1}{2(n-1)}}.
$$
\end{lemma}
\begin{proof}
The proof of this lemma is almost identical to the one in \cite[Theorem 22.8]{Ma}. Indeed, it follows from the density lower bound (see Theorem \ref{Density lower bound}), the relative isoperimetric inequality and the compactness result proved in the previous section.
\end{proof}
Proceeding as in \cite{Ma}, we give the following Lipschitz approximation lemma, which is a consequence of the height bound lemma. Its proof follows exactly as in \cite[Theorem 23.7]{Ma}. {It is a foundamental step in the long journey to the regularity because it provides a connection between the regularity theories for parametric and non-parametric variational problems. Indeed we are able to prove for $(\Lambda,\alpha)$-minimizers that the smallness of the excess guaranties that $\partial E$ can be locally almost entirely covered by the graph of a Lipschitz function}.
\begin{theorem}[Lipschitz approximation]\label{LipApp}
Let $(E,u)$ be a $(\Lambda,\alpha)$-minimizer of $\mathcal{F}$ in $B_r(x_0)$.
There exist two positive constants $\varepsilon_3$ and $C_8$, depending on $\norm{\D u}_{L^2(B_r(x_0))}$, such that if $x_0\in \partial E$ and
$$
{\mathbf e}(x_0,r,e_n)<\varepsilon _3,
$$
then there exists a Lipschitz function $f:\R^{n-1}\rightarrow \R$ such that
$$
\sup_{x'\in \R^{n-1}}\frac{|f(x')|}{r}\leq C_8{\mathbf e}(x_0,r,e_n)^{\frac{1}{2(n-1)}},\quad \norm{\nabla'f}_{L^{\infty}}\leq 1,
$$
and
$$
\frac{1}{r^{n-1}}\mathcal{H}^{n-1}((\partial E \Delta \Gamma_f)\cap B_{r/2}(x_0))\leq C_8 {\mathbf e}(x_0,r,e_n),
$$
where $\Gamma_f$ is the graph of $f$. Moreover,
$$
\frac{1}{r^{n-1}}\int_{D_{r/2}(x'_0)}|\nabla'f|^2\,dx'\leq C_8 {\mathbf e}(x_0,r,e_n).
$$
\end{theorem}

\section{Reverse Poincaré inequality}
\label{Reverse Poincaré inequality}
{In this section we shall prove a reverse Poincaré inequality. This is the counterpart for $(\Lambda,\alpha)$-minimizers of the well-known Caccioppoli inequality for weak solutions of elliptic equations. The proof of the results of this section can be obtained as in the case of $\Lambda$-minimizers of the perimeter (sse \cite[Section 24]{Ma}). For the sake of completeness we give here the main steps of the proof underlining the minor changes.} We will need first a weak form.
\begin{lemma}[Weak reverse Poincaré inequality]
\label{Weak Reverse Poincaré}
If $(E,u)$ is a $(\Lambda,\alpha)$-minimizer of $\mathcal{F}$ in $\mathbf{C}_4$ such that
\begin{equation*}
|x_n|<\frac{1}{8},\quad\forall x\in \mathbf{C}_2\cap\dd E,
\end{equation*}
\begin{equation*}
\left|\left\{ x\in \mathbf{C}_2\setminus E\,:\, x_n<-\frac{1}{8} \right\}\right|=\left|\left\{ x\in \mathbf{C}_2\cap E\,:\, x_n>\frac{1}{8} \right\}\right|=0,
\end{equation*}
and if $z\in\R^{n-1}$ and $s>0$ are such that
\begin{equation}
\label{eqqq18}
\mathbf{K}_s(z)\subset \mathbf{C}_2,\qquad \mathcal{H}^{n-1}(\dd E\cap\dd \mathbf{K}_s(z))=0,
\end{equation}
then, for every $|c|<\frac{1}{4}$,
\begin{align*}
& P(E;\mathbf{K}_{\frac{s}{2}}(z))-\mathcal{H}^{n-1}(D_{\frac{s}{2}}(z))\leq C(n,N,L)\Bigg\{\bigg[\left(P(E;\mathbf{K}_s(z))-\mathcal{H}^{n-1}(D_s(z))\right)\\
& \times\int_{\mathbf{K}_s(z)\cap\dd^*E}\frac{(x_n-c)^2}{s^2}d\,\mathcal{H}^{n-1}\bigg]^{\frac{1}{2}}+\Lambda s^{(n-1)\alpha}+\int_{\mathbf{K}_{s}}|\D u|^2\,dx\Bigg\}.
\end{align*}
\end{lemma}
\begin{proof}
We may assume $z=0$.\\
\textbf{Step 1:} The set function
\begin{equation*}
\zeta(G)=P(E;\mathbf{C}_2\cap p^{-1}(G))-\mathcal{H}^{n-1}(G),\quad\text{for }G\subset D_2,
\end{equation*}
defines a Radon measure on $\R^{n-1}$, concentrated on $D_2$.\\
\textbf{Step 2: } Since $E$ is a set of locally finite perimeter, by \cite[Theorem 13.8]{Ma} there exist a sequence $\{E_h\}_{h\in\N}$ of open subsets of $\R^n$ with smooth boundary and a vanishing sequence $\{\varepsilon_h\}_{h\in\N}\subset\R^+$ such that
\begin{equation*}
E_h\overset{loc}{\rightarrow}E,\quad\mathcal{H}^{n-1}\llcorner\dd E_h\rightarrow\mathcal{H}^{n-1}\llcorner\dd E,\quad \dd E_h\subset I_{\varepsilon_h}(\dd E),
\end{equation*}
as $h\rightarrow+\infty$, where $I_{\varepsilon_h}(\dd E)$ is a tubolar neighborhood of $\dd E$ with half-lenght $\varepsilon_h$. By Coarea formula we get
\begin{equation*}
\mathcal{H}^{n-1}(\dd \mathbf{K}_{rs}\cap(E^{(1)}\Delta E_h))\rightarrow0,\quad \text{for a.e. }r\in\left(\frac{2}{3},\frac{3}{4}\right).
\end{equation*}
Moreover, provided $h$ is large enough, by $\dd E_h\subset I_{\varepsilon_h}(\dd E)$, we get:
\begin{equation*}
|x_n|<\frac{1}{4},\quad\forall x\in \mathbf{C}_2\cap\dd E_h,
\end{equation*}
\begin{equation*}
\left\{ x\in \mathbf{C}_2\,:\, x_n<-\frac{1}{4} \right\}\subset \mathbf{C}_2\cap E_h\subset\left\{ x\in \mathbf{C}_2\,:\, x_n<\frac{1}{4} \right\}.
\end{equation*}
Therefore, given $\lambda\in\left(0,\frac{1}{4}\right)$ and $|c|<\frac{1}{4}$, we are in position to apply \cite[Lemma 24.8]{Ma} to every $E_h$ to deduce that there exists $I_h\subset\left( \frac{2}{3},\frac{3}{4} \right)$, with $|I_h|\geq\frac{1}{24}$, and, for any $r\in I_h$, there exists an open subset $F_h$ of $\R^n$ of locally finite perimeter such that
\begin{equation}
\label{eqqq15}
F_h\cap\dd \mathbf{K}_{rs}=E_h\cap\dd \mathbf{K}_{rs},\\
\end{equation} 
\begin{equation}
\label{eqqq16}
\mathbf{K}_{\frac{r}{2}}\cap\dd F_h=D_{\frac{s}{2}}\times\{c\},
\end{equation}
\begin{align}
\label{eqqq17}
P(F_h;\mathbf{K}_{rs})-\mathcal{H}^{n-1}(D_{rs})\leq 
& c(n)\bigg\{ \lambda\left( P(E_h;\mathbf{K}_s)-\mathcal{H}^{n-1}(D_s) \right)+\frac{1}{\lambda}\int_{\mathbf{K}_s\cap\dd E_h}\frac{|x_n-c|^2}{s^2}\,d\mathcal{H}^{n-1} \bigg\}.
\end{align}
Clearly $\displaystyle\bigcap_{h\in\N}\bigcup_{k\geq h}|I_k|\geq\frac{1}{24}>0$ and thus there exist a divergent subsequence $\{h_k\}_{k\in\N}$ and $r\in\left( \frac{2}{3},\frac{3}{4}\right)$ such that
\begin{equation*}
r\in\bigcap_{k\in\N} I_{h_k} \quad\text{and}\quad \lim_{k\rightarrow+\infty}\mathcal{H}^{n-1}(\dd \mathbf{K}_{rs}\cap(E^{(1)}\Delta E_{h_k}))=0.
\end{equation*}
We will write $F_k$ in place of $F_{h_k}$. Now we test the $(\Lambda,\alpha)$-minimality of $(E,u)$ in $\mathbf{C}_4$ with $(G_k,u)$, where $G_k=(F_k\cap \mathbf{K}_{rs})\cup (E\setminus \mathbf{K}_{rs})$, as $E\Delta G_k\subset\subset \mathbf{K}_s\subset\subset B_4$. By \cite[(16.33)]{Ma} we infer:
\begin{align*}
P(E;\mathbf{K}_{rs})
& \leq P(G_k;\mathbf{K}_{rs})+\Lambda|(E\Delta F_k)\cap \mathbf{K}_{rs}|^\alpha+\int_{\mathbf{K}_{rs}}G(x,u,\D u)[\mathbbm{1}_{G_k}-\mathbbm{1}_E]\,dx\\
& \leq P(F_k;\mathbf{K}_{rs})+\sigma_k+\Lambda|(E\Delta F_k)\cap \mathbf{K}_{rs}|^\alpha+c(n,N,L)\int_{\mathbf{K}_{rs}}(|\D u|^2+1)\,dx,
\end{align*}
with $\sigma_k=\mathcal{H}^{n-1}(\dd \mathbf{K}_{rs}\cap (E^{(1)}\Delta F_k))=\mathcal{H}^{n-1}(\dd \mathbf{K}_{rs}\cap (E^{(1)}\Delta E_{h_k}))\rightarrow0$, thanks to \eqref{eqqq15}, as $k\rightarrow+\infty$. Thus, since $\zeta$ is nondecreasing and $r\geq\frac{2}{3}$, by \eqref{eqqq17} we deduce that
\begin{align*}
& P(E;\mathbf{K}_{\frac{s}{2}})-\mathcal{H}^{n-1}(D_{\frac{s}{2}})=\zeta(D_{\frac{s}{2}})\leq\zeta(D_{rs})=P(E;\mathbf{K}_{rs})-\mathcal{H}^{n-1}(D_{rs})\\
& \leq P(F_k;\mathbf{K}_{rs})-\mathcal{H}^{n-1}(D_{rs})+\sigma_k+\Lambda|(E\Delta F_k)\cap \mathbf{K}_{rs}|^\alpha+c(n,N,L)\int_{\mathbf{K}_{rs}}(|\D u|^2+1)\,dx\\
& \leq c(n)\bigg\{ \lambda\left( P(E_{h_k};\mathbf{K}_s)-\mathcal{H}^{n-1}(D_s) \right)+\frac{1}{\lambda}\int_{\mathbf{K}_s\cap\dd E_{h_k}}\frac{|x_n-c|^2}{s^2}\,d\mathcal{H}^{n-1} \bigg\}\\
& +c(n,N,L)\left(\Lambda s^{(n-1)\alpha}+\int_{\mathbf{K}_{s}}|\D u|^2\,dx\right).
\end{align*}
Letting $k\rightarrow +\infty$, \eqref{eqqq18} implies that $P(E_{h(k)};\mathbf{K}_s)\rightarrow P(E;\mathbf{K}_s)$ and therefore
\begin{align}
\label{eqqq19}
P(E;\mathbf{K}_{\frac{s}{2}})-\mathcal{H}^{n-1}(D_{\frac{s}{2}})\notag
& \leq c(n)\bigg\{ \lambda\left( P(E;\mathbf{K}_s)-\mathcal{H}^{n-1}(D_s) \right)+\frac{1}{\lambda}\int_{\mathbf{K}_s\cap\dd E}\frac{|x_n-c|^2}{s^2}\,d\mathcal{H}^{n-1} \bigg\}\notag\\
& +c(n,N,L)\left(\Lambda s^{(n-1)\alpha}+\int_{\mathbf{K}_{s}}|\D u|^2\,dx\right),
\end{align}
for any $\lambda\in\left(0,\frac{1}{4}\right)$. If $\lambda>\frac{1}{4}$,
\begin{align*}
& P(E;\mathbf{K}_{\frac{s}{2}})-\mathcal{H}^{n-1}(D_{\frac{s}{2}})=\zeta(D_{\frac{s}{2}})\leq\zeta(D_{rs})\\
& \leq 4\lambda P(E;\mathbf{K}_{rs})-\mathcal{H}^{n-1}(D_{rs})\leq c(n)\lambda\left( P(E;\mathbf{K}_s)-\mathcal{H}^{n-1}(D_s) \right)
\end{align*}
and thus \eqref{eqqq19} holds  true for $\lambda>0$, provided we choose $c(n)\geq 4$. Minimizing over $\lambda$, we get the thesis.
\end{proof}

\begin{theorem}[Reverse Poincaré Inequality] 
There exists a positive constant $C_9=C_9(n,N,L,\alpha)$ such that if $(E,u)$ be a $(\Lambda,\alpha)$-minimizer of $\mathcal{F}$ in $\mathbf{C}_{4r}(x_0,\nu)$ with $x_0\in \partial E$ and
\begin{equation*}
\ecc(x_0,4r,\nu)<\omega\left(n,\frac{1}{8}\right),
\end{equation*}
then
\begin{equation*}
\ecc(x_0,r,\nu)\leq
C_9\biggl(\frac{1}{r^{n+1}}\int_{\partial E\cap \mathbf{C}_{2r}(x_0,\nu)}|\left<\nu,x-x_0\right>-c|^2d\mathcal{H}^{n-1}+\Lambda r^{\gamma}+\frac{1}{r^{n-1}}\int_{\mathbf{K}_{2r}}|\D u|^2\,dx\biggr),
\end{equation*}
for every $c\in \mathbbm R$.
\end{theorem}
\begin{proof}
Up to replacing $(E,u)$ with {$\left(\frac{E-x_0}{r},r^{-\frac{1}{2}}u(x_0+ry)\right)$, see Lemma \ref{Lemma riscalamento}}, we can assume that $(E,u)$ is a $(\Lambda r^{\gamma},\alpha)$-minimizer of $\mathcal{F}_r$ in $\mathbf{C}_4$, $0\in\dd E$ and, by \cite[Proposition 22.1]{Ma},
\begin{equation*}
\ecc(0,4,e_n)\leq \omega\left(n,\frac{1}{8}\right).
\end{equation*}
Applying \cite[Lemma 22.10 and Lemma 22.11]{Ma}, we get that
\begin{equation*}
|x_n|<\frac{1}{4},\quad\forall x\in \mathbf{C}_2\cap\dd E,
\end{equation*}
\begin{equation*}
\left|\left\{ x\in \mathbf{C}_2\setminus E\,:\, x_n<-\frac{1}{8} \right\}\right|=\left|\left\{ x\in \mathbf{C}_2\cap E\,:\, x_n>\frac{1}{8} \right\}\right|=0.
\end{equation*}
\begin{equation*}
\mathcal{H}^{n-1}(G)=\int_{\mathbf{C}_2\cap\dd^*E\cap p^{-1}(G)}\langle\nu_E, e_n\rangle\,d\mathcal{H}^{n-1},\quad\forall G\subset D_2.
\end{equation*}
Since
\begin{align*}
\textbf{e}_n(1)
& =\int_{\mathbf{C}_1\cap\dd^* E}(1-\langle\nu_E, e_n\rangle)\,d\mathcal{H}^{n-1}=P(E;\mathbf{C}_1)-\int_{\mathbf{C}_1\cap\dd^* E}\langle\nu_E, e_n\rangle\,d\mathcal{H}^{n-1}=P(E;\mathbf{C}_1)-\mathcal{H}^{n-1}(D_1),
\end{align*}
then our aim is to show
\begin{equation*}
P(E;\mathbf{C}_1)-\mathcal{H}^{n-1}(D_1)\leq C_{9}\bigg( \int_{\mathbf{C}_2\cap\dd E}|x_n-c|^2\,d\mathcal{H}^{n-1}+\Lambda r^{\gamma}+\int_{\mathbf{K}_{2}}|\D u|^2\,dx \bigg),
\end{equation*}
for any $c\in\R$. Actually it suffices to prove it only for $|c|<\frac{1}{4}$; indeed, for $|c|\geq\frac{1}{4}$, we have:
\begin{equation}
\int_{\mathbf{C}_2\cap\dd E}|x_n-c|^2\,d\mathcal{H}^{n-1}\geq \int_{\mathbf{C}_2\cap\dd E}(|c|-|x_n|)^2\,d\mathcal{H}^{n-1}\geq \frac{1}{64}P(E;\mathbf{C}_2)\geq \frac{1}{64}P(E;\mathbf{C}_1).
\end{equation}
\textbf{Step 2:} the set function $\zeta(G)=P(E;\mathbf{C}_2\cap p^{-1}(G))-\mathcal{H}^{n-1}(G)$, for $G\subset D_2$, defines a Radon measure on $\R^{n-1}$, concentrated on $D_2$. We apply Lemma \ref{Weak Reverse Poincaré} to $E$ in every cylinder $\mathbf{K}_s(z)$ with $z\in\R^{n-1}$ and $s>0$ such that
\begin{equation}
\label{eqqq20}
D_{2s}(z)\subset D_2,\qquad \mathcal{H}^{n-1}(\dd E\cap\dd \mathbf{K}_{2s}(z))=0,
\end{equation}
to get that
\begin{equation*}
\zeta(D_s(z))\leq C(n,N,L,\alpha)\left\{ (\zeta(D_{2s}(z))h)^{\frac{1}{2}}
+\Lambda r^{\gamma}s^{(n-1)\alpha}+\int_{\mathbf{K}_{2s}(z)}|\D u|^2\,dx \right\},
\end{equation*}
where
\begin{equation*}
h:=\inf_{|c|<\frac{1}{4}}\int_{\mathbf{C}_2\cap\dd E}|x_n-c|^2\,d\mathcal{H}^{n-1}.
\end{equation*}
Multiplying by $s^2$ and using an approximation argument to remove the second assumption in \eqref{eqqq20}, we obtain:
\begin{equation}
\label{eqqq21}
s^2\zeta(D_s(z))\leq c(n,N,L,\alpha)\left(\sqrt{s^2\zeta(D_{2s}(z))h}+\Lambda r^{\gamma}+\int_{\mathbf{K}_{2s}(z)}|\D u|^2\,dx\right), 
\end{equation}
for $D_{2s}(z)\subset D_2$, where we used that $s<1$. In order to prove the thesis, we use a covering argument by setting
\begin{equation*}
Q=\sup_{D_{2s}(z)\subset D_2}s^2\zeta(D_s(z))< +\infty.
\end{equation*}
We cover $D_s(z)$ by finitely many balls $\{ D\left(z_k,\frac{s}{4}\right)\}_{k\in\{1,\dots,\tilde{N}\}}$ with centers $z_k\in D_s(z)$. Of course, this can be done with $\tilde{N}\leq \tilde{N}(n)$, for some $\tilde{N}(n)\in\N$. Hence, by the sub-additivity of $\zeta$ and \eqref{eqqq21} for $\frac{s}{4}$, since $D_{s}(z_k)\subset D_2$, we have:
\begin{align*}
s^2\zeta(D_s(z))
& \leq s^2\sum_{k=1}^{\tilde{N}} \zeta\left(D_{\frac{s}{4}}(z_k)\right)=16\sum_{k=1}^{\tilde{N}}\left(\frac{s}{4}\right)^2\zeta\left(D_{\frac{s}{4}}(z_k)\right)\\
& \leq c(n,N,L,\alpha)\sum_{k=1}^{\tilde{N}}\Bigg( \sqrt{\left( \frac{s}{2}\right)^2\zeta\left( D_{\frac{s}{2}}(z_k)\right)h}+\Lambda r^{\gamma}+\int_{\mathbf{K}_{2s}(z)}|\D u|^2\,dx \Bigg)\\
& \leq c(n,N,L,\alpha)\left(\sqrt{Qh}+\Lambda r^{\gamma}+\int_{\mathbf{K}_{2s}(z)}|\D u|^2\,dx\right).
\end{align*}
Passing to the supremum for $D_{2s}(z)\subset D_2$ we infer that
\begin{equation*}
Q\leq c(n,N,L,\alpha)\left(\sqrt{Qh}+\Lambda r^{\gamma}+\int_{\mathbf{K}_2}|\D u|^2\,dx\right).
\end{equation*}
If $\sqrt{Qh}\leq\Lambda r^{\gamma}+\int_{\mathbf{K}_2}|\D u|^2\,dx$, then $Q\leq c(n,N,L,\alpha)\left(\Lambda r^{\gamma}+\int_{\mathbf{K}_2}|\D u|^2\,dx\right)$.\\ If $\sqrt{Qh}>\Lambda r^{\gamma}+{\int_{\mathbf{K}_2}|\D u|^2\,dx}$, then $Q\leq c(n,N,L,\alpha)\sqrt{Qh}$ and thus $Q\leq c(n,N,L,\alpha)h$.\\ In both cases we obtain:
\begin{equation*}
Q\leq c(n,N,L,\alpha)\left(h+\Lambda r^{\gamma}+\int_{\mathbf{K}_2}|\D u|^2\,dx\right),
\end{equation*}
which leads to the thesis.
\end{proof}

%

\section{Energy first variation}
\label{Energy first variation}
In this section we deduce a kind of Taylor's expansion formula, with respect to a parameter $t\in \R$, for the energy quantity involved in the definition of $(\Lambda,\alpha)$-minimizer, under a ``small'' domain perturbation of the type
$$
\Phi_t(x)=x+tX(x).
$$
We start with the energy of the rescaled functional $\mathcal{F}_r$. For the sake of simplicity we will denote with $A_1(x,s)$ the matrix whose entries are $a_{hk}(x,s)$, $A_2(x,s)$ the vector of components $a_{h}(x,s)$, $A_3(x,s)=a(x,s)$ and similarly for $B_i$, $i=1,2,3$.
Then we define
\begin{equation*}
\begin{split}
&\mathcal{F}_r(w;D)
:=\int_{B_1}\big[ F_r(x,w,\D w)+\mathbbm{1}_DG_r(x,w,\D w) \big]\,dx\\
& = \int_{B_1}\big[ \langle (A_{1r}+\mathbbm{1}_D B_{1r})\D w,\D w \rangle + \sqrt{r}\langle A_{2r}+\mathbbm{1}_D B_{2r},\D w \rangle + r(A_{3r}+\mathbbm{1}_DB_{3r}) \big]\,dx,
\end{split}
\end{equation*}
where $r>0$, $x_0\in\Omega$,  $A_{ir}(y,w):=A_i (x_0+ry,\sqrt{r} w)$, $B_{ir}(y,w):=B_i (x_0+ry,\sqrt{r} w)$, for $i=1,2,3$.

\begin{theorem}[First variation of the bulk term]\label{VarB}
Let $u\in H^1(B_1)$ and let us fix
$X\in C_0^1(B_1;\R^n)$. We define $\Phi_t(x):=x+tX(x)$, for any $t>0$. Accordingly we define
$$E_t:=\Phi_t(E),\quad u_t:=u\circ\Phi_t^{-1}.$$ There exists a constant 
$\overline{C}=\overline{C}(N,L,L_{\alpha},\norm{X}_{\infty},\norm{\D X}_{\infty})>0$ such that

\begin{align}
\label{eqq1}
& \int_{B_1}\big[F_r(y,u_t,\D u_t)+\mathbbm{1}_{E_t}(y)G_r(y,u_t,\D u_t)\big]\,dy-\int_{B_1}\big[F_r(x,u,\D u)+\mathbbm{1}_{E}(x)G_r(x,u,\D u)\big]\,dx\notag\\
& \leq \overline{C}(t^{\alpha}+o(t))\int_{B_1}\big( |\D u|^2+r \big)\,dx,
\end{align}
where $L_{\alpha}$ is defined in \eqref{Hoelderianity2}.
\end{theorem}

\begin{proof}
Taking into account that
\begin{equation*}
\D \Phi_t^{-1}(\Phi_t(x))=I-t\D X(x)+o(t), \quad \textnormal{J}\Phi_t(x)=1+t\textnormal{div}X(x)+o(t).
\end{equation*}
we obtain:
\begin{align*}
& \int_{B_1}\big[F_r(y,u_t,\D u_t)+\mathbbm{1}_{E_t}(y)G_r(y,u_t,\D u_t)\big]\,dy\\
& =\int_{B_1} \big[ F_r(\Phi_t(x),u,\D u)+\mathbbm{1}_E(x) G_r(\Phi_t(x),u,\D u) \big](1+t\textnormal{div}X+o(t))\,dx\\
& - (t+o(t))\int_{B_1}\big[  t\big\langle C_1\D u\D X,\D u\D X \big\rangle +  2\big\langle C_1\D u\D X,\D u \big\rangle +\sqrt{r}\big\langle C_2,\D u\D X \big\rangle\big]\,dx,
\end{align*}
where we set
\begin{equation*}
C_i:=\tilde{A}_{ir}+\mathbbm{1}_{E}\tilde{B}_{ir}=A_{ir}(\Phi_t(x),u)+\mathbbm{1}_{E}(x)B_{ir}(\Phi_t(x),u),
\end{equation*}
for $i=1,2,3$. From the previous identity, by subtracting  the term
$$\int_{B_1}\big[F_r(x,u,\D u)+\mathbbm{1}_{E}(x)G_r(x,u,\D u)\big]\,dx,$$
we gain:
\begin{align}
\label{eqq1}
& \int_{B_1}\big[F_r(y,u_t,\D u_t)+\mathbbm{1}_{E_t}(y)G_r(y,u_t,\D u_t)\big]\,dy 
 -\int_{B_1}\big[F_r(x,u,\D u)+\mathbbm{1}_{E}(x)G_r(x,u,\D u)\big]\,dx\notag\\
& = \bigg[\int_{B_1} \big[ F_r(\Phi_t(x),u,\D u)+\mathbbm{1}_E(x) G_r(\Phi_t(x),u,\D u)-[F_r(x,u,\D u)+\mathbbm{1}_{E}(x)G_r(x,u,\D u)] \big]\,dx\bigg]\notag \\
& + \bigg[t\int_{B_1} \big[F_r(\Phi_t(x),u,\D u)+\mathbbm{1}_E(x) G_r(\Phi_t(x),u,\D u)\big]\textnormal{div}X\,dx \notag\\
&+o(t)\int_{B_1} \big[F_r(\Phi_t(x),u,\D u)+\mathbbm{1}_E(x) G_r(\Phi_t(x),u,\D u)\big]\textnormal\,dx \notag\\
& - (t+o(t))\int_{B_1}\big[ t\big\langle C_1\D u\D X,\D u\D X \big\rangle + 2\big\langle C_1\D u\D X,\D u \big\rangle + \sqrt{r}\big\langle C_2,\D u\D X \big\rangle]\,dx\bigg]=:\big[I_1\big]+\big[I_2\big]\notag.
\end{align}
Let us estimate separately the two terms $I_1,I_2$ on the right-hand side. By the H\"{o}lder continuity of the data with respect to the first variable given in \eqref{Hoelderianity2} and Young's inequality we get
\begin{equation*}
\begin{split}
I_1
& = \int_{B_1} \big[ \big\langle F_r(\Phi_t(x),u,\D u)+\mathbbm{1}_E(x) G_r(\Phi_t(x),u,\D u)-[F_r(x,u,\D u)+\mathbbm{1}_{E}(x)G_r(x,u,\D u)] \big]\,dx\notag\\
& \leq c(L_{\alpha})t^{\alpha}\int_{B_1}|X|[|\D u|^2+\sqrt{r}|\D u|+r]\,dx \leq c(L_{\alpha},\norm{X}_\infty)t^{\alpha}\int_{B_1}[|\D u|^2+r]\,dx.
\end{split}
\end{equation*}
Regarding $I_2$ we have that
\begin{align}
\label{eqq4}
I_2 &\leq  (t+o(t))(1+\norm{\nabla X}_{\infty})\int_{B_1} \big|F_r(\Phi_t(x),u,\D u)+\mathbbm{1}_E(x) G_r(\Phi_t(x),u,\D u)\big|\,dx\notag\\
& +(t+o(t))(1+\norm{\nabla X}_{\infty})^2\int_{B_1}\big| t\big\langle C_1\D u,\D u \big\rangle + 2\big\langle C_1\D u,\D u \big\rangle + \sqrt{r}\big\langle C_2,\D u\big\rangle\big|\,dx\notag\\
&\leq {C}(t+o(t))\int_{B_1}\big( |\D u|^2+r \big)\,dx,
\end{align}
where ${C}={C}(N,L,\norm{\D X}_{\infty})$.
From the last estimates the thesis easly follows.
\end{proof}

The second estimate concerns the perimeter (see \cite[Theorem 17.5]{Ma}).
\begin{theorem}[First variation of the perimeter]\label{VarP}
If $A\subset\R^n$ is an open set, $E\subset\R^n$ is a set of locally finite perimeter and $\Phi_t(x):=x+tX(x)$ for some fixed  $X\in C_0^1(A;\R^n)$, then
\begin{equation}
P(\Phi_t(E);A)-P(E;A)=(t+O(t^2))\int_{\dd^*E} \textnormal{div}_E X\,d\mathcal{H}^{n-1},
\end{equation}
where the tangential divergence of $X$, $\textnormal{div}_E X:\partial^* E\rightarrow \R$, is the Borel function defined as
\begin{equation}
\textnormal{div}_E X=\textnormal{div} X -\langle\nu_{E}, \nabla X\nu_E\rangle.
\end{equation}
\end{theorem}
The last result we will use in the sequel concerns the penalization term (see \cite[Lemma 17.9]{Ma}).
\begin{theorem}\label{VarD}
Let $A\subset\R^n$ be an open set, $E\subset\R^n$ be a set of locally finite perimeter and $\Phi_t(x):=x+tX(x)$, for some fixed $X\in C_0^1(A;\R^n)$, be a local variation in $A$, i.e. $\left\{x\neq \Phi_t(x)\right\}\subset K\subset A$, for some compact set $K\subset A$ and for $|t|<\varepsilon_0$. Then
\begin{equation}
|\Phi_t(E)\Delta E|\leq C|t|P(E;K),
\end{equation}
where $C$ is a positive constant.
\end{theorem}
\section{Excess improvement}
\label{Excess improvement}

\begin{theorem}[Excess improvement]
\label{Miglioramento eccesso}
For every $\tau\in \big(0,\frac{1}{2}\big)$ and $M>0$ there exists a constant $\varepsilon_4=\varepsilon_4(\tau,M)\in(0,1)$ such that if $(E,u)$ is a $(\Lambda,\alpha)$-minimizer of $\mathcal{F}$ in $B_r(x_0)$ with $x_0\in\dd E$ and
\begin{equation}
\label{3}
{\mathbf e}(x_0,r)\leq\varepsilon_4, \quad\mathcal{D}(x_0,r)+r^\gamma\leq M{\mathbf e}(x_0,r),
\end{equation}
then there exists a positive constant $C_{10}$, depending on $\norm{\D u}_{L^2(B_r(x_0))}$, such that
\begin{equation*}
{\mathbf e}(x_0,\tau r)\leq C_{10}(\tau^2{\mathbf e}(x_0,r)+\mathcal{D}(x_0,2\tau r)+(\tau r)^\gamma).
\end{equation*}
\end{theorem}
\begin{proof}
Without loss of generality we may assume that $\tau<\frac{1}{8}$. Let us rescale 
and assume by contradiction that there exist an infinitesimal sequence $\lbrace\varepsilon_h\rbrace_{h\in\N}\subseteq\R^+$, a sequence $\lbrace r_h\rbrace_{h\in\N}\subseteq\R^+$ and a sequence $\lbrace (E_h,u_h) \rbrace_{h\in\N}$ of $(\Lambda r_h^{\gamma},\alpha)$-minimizers of $\mathcal{F}_{r_h}$ in $B_1$, with equibounded energies, such that, denoting by ${\mathbf e}_h$ the excess of $E_h$ and by $\mathcal{D}_h$  the rescaled Dirichlet integral of $u_h$, we have
\begin{equation*}
{\mathbf e}_h(0,1)=\varepsilon_h, \quad \mathcal{D}_h(0,1)+r_h^{\gamma}\leq M\varepsilon_h
\end{equation*}
and

\begin{equation*}
{\mathbf e}_h(0,\tau)>C_{10}(\tau^2{{\mathbf e}_h(0,1)}+{\mathcal{D}_h}(0,2\tau)+(\tau r_h)^{\gamma}),
\end{equation*}
with some positive constant $C_{10}$ to be chosen.
Up to rotating each $E_h$ we may also assume that, for all $h\in\N$,
\begin{equation*}
{\mathbf e}_h(0,1)=\frac{1}{2}\int_{\dd E_h\cap B_1}\abs{\nu_{E_h}-e_n}^2\, d\mathcal{H}^{n-1}.
\end{equation*}
\textbf{Step 1.} Thanks to the Lipschitz approximation theorem, for $h$ sufficiently large, there exists a 1-Lipschitz function $f_h\colon\R^{n-1}\rightarrow\R$ such that
\begin{equation}
\label{1}
\sup_{\R^{n-1}}\abs{f_h}\leq C_8\varepsilon_h^{\frac{1}{2(n-1)}}, \quad \mathcal{H}^{n-1}((\dd E_h\Delta\Gamma_{f_h})\cap B_{\frac{1}{2}})\leq C_8\varepsilon_h, \quad \int_{ {D_{\frac{1}{2}}}}\abs{\D' f_h}^2\,dx'\leq C_8\varepsilon_h.
\end{equation}
We define
\begin{equation*}
g_h(x'):=\frac{f_h(x')-a_h}{\sqrt{\varepsilon_h}}, \quad\text{where}\quad a_h=\fint_{ {D_{\frac{1}{2}}}}f_h\,dx'
\end{equation*}
and we assume, up to a subsequence, that $\lbrace g_h \rbrace_{h\in\N}$ converges weakly in $H^1( {D_{\frac{1}{2}}})$ and strongly in $L^2( {D_{\frac{1}{2}}})$ to a function $g$.\\
We prove that $g$ is harmonic in $ {D_{\frac{1}{2}}}$. It is enough to show that
\begin{equation}
\label{2}
\lim_{h\rightarrow +\infty}\frac{1}{\sqrt{\varepsilon_h}}\int_{ {D_{\frac{1}{2}}}}\frac{\langle \D' f_h,\D'\phi\rangle}{\sqrt{1+\abs{\D' f_h}^2}}\,dx'=0,
\end{equation}
for all $\phi\in C_0^1( {D_{\frac{1}{2}}})$; indeed, if $\phi\in C_0^1( {D_{\frac{1}{2}}})$, by weak convergence we have
\begin{equation*}
\begin{split}
& \int_{ {D_{\frac{1}{2}}}} \langle\D' g,\D'\phi\rangle\,dx' =
\lim_{h\rightarrow +\infty}\frac{1}{\sqrt{\varepsilon_h}}\int_{ {D_{\frac{1}{2}}}}\langle\D' f_h,\D'\phi\rangle\,dx'\\
& = \lim_{h\rightarrow +\infty}\frac{1}{\sqrt{\varepsilon_h}}\bigg\lbrace \int_{ {D_{\frac{1}{2}}}}\frac{\langle \D 'f_h,\D'\phi\rangle}{\sqrt{1+\abs{\D' f_h}^2}}\,dx'+\int_{ {D_{\frac{1}{2}}}} \bigg[ \langle \D' f_h,\D'\phi\rangle - \frac{\langle \D' f_h,\D'\phi\rangle}{\sqrt{1+\abs{\D' f_h}^2}}\bigg]\,dx' \bigg\rbrace.
\end{split}
\end{equation*}
Using the Lipschitz continuity of $f_h$ and the third inequality in \eqref{1}, we infer that the second term in the previous equality is infinitesimal:
\begin{equation*}
\begin{split}
& \limsup_{h\rightarrow +\infty}\frac{1}{\sqrt{\varepsilon_h}}\bigg |\int_{ {D_{\frac{1}{2}}}} \bigg [ \langle \D' f_h,\D'\phi\rangle - \frac{\langle \D' f_h,\D'\phi\rangle}{\sqrt{1+\abs{\D' f_h}^2}}\bigg]\,dx'\bigg |\\
& \leq \limsup_{h\rightarrow +\infty}\frac{1}{\sqrt{\varepsilon_h}} \int_{ {D_{\frac{1}{2}}}} \abs{\D' f_h}\abs{\D'\phi}\frac{\sqrt{1+\abs{\D' f_h}^2}-1}{\sqrt{1+\abs{\D' f_h}^2}}\,dx'\\
& \leq \limsup_{h\rightarrow +\infty}\frac{1}{\sqrt{\varepsilon_h}} \int_{ {D_{\frac{1}{2}}}} \abs{\D'\phi}\abs{\D' f_h}^2\,dx'\leq \lim_{h\rightarrow +\infty}C_8\norm{\D'\phi}_\infty\sqrt{\varepsilon_h}=0.
\end{split}
\end{equation*}
Therefore, we should prove \eqref{2}. We fix $\delta>0$ so that spt$\,\phi\times [-2\delta,2\delta]\subset B_{\frac{1}{2}}$ and choose a cut-off function $\psi\colon\R\rightarrow[0,1]$ with spt$\,\psi\subset (-2\delta,2\delta)$, $\psi=1$ in $(-\delta,\delta)$. Let us define
$$
\Phi_{\varepsilon_h}(x):=x+\varepsilon_h X(x),\quad \text{where }X(x)=\phi(x')\psi(x_n) e_n,
$$
and 
$$\tilde{E}_h:=\Phi_{\varepsilon_h}(E_h), \quad\tilde{u}_h:=u\circ \Phi^{-1}_{\varepsilon_h}.$$
By the $(\Lambda,\alpha)$-minimality of $(E_h,u_h)$ we deduce that
$$
\mathcal{F}_{r_h}(E_h,u_h)\leq \mathcal{F}_{r_h}(\tilde{E}_h,\tilde{u}_h)+\Lambda r_h^{\gamma}|\tilde{E}_h\Delta E_h|^{\alpha}.
$$
Then we may estimate
\begin{align}
\label{Minuh}
& P(E_h;B_{\frac 12})-P(\tilde{E}_h;B_{\frac 12})\notag \\
&\leq \int_{B_\frac{1}{2}}\big[F_r(y,\tilde{u}_h,\D \tilde{u}_h)+\mathbbm{1}_{\tilde{E}_h}(y)G_r(y,\tilde{u}_h,\D \tilde{u}_h)\big]\,dy-\int_{B_\frac{1}{2}}\big[F_r(x,u,\D u)+\mathbbm{1}_{E}(x)G_r(x,u,\D u)\big]\,dx\\
& + \Lambda r_h^{\gamma}|\Phi_{\varepsilon_h}(E_h)\Delta E_h|^{\alpha}.\\
\end{align}
Applying Theorem \ref{VarB} and Theorem \ref{VarD} in the right-hand side we get
\begin{align}
\label{Ar}
& P(E_h;B_{\frac 12})-P(\tilde{E}_h;B_{\frac 12})\leq C\bigg[\bigl(\varepsilon_h^{\alpha}+o(\varepsilon_h)\bigr)\int_{B_1}\big( |\D u_h|^2+r_h \big)\,dx
+r_h^{\gamma}\varepsilon_h^{\alpha}(P(E_h;B_1))^{\alpha}\bigg],
\end{align}
for some $C=C(N,L,L_{\alpha},\alpha,\Lambda,\norm{X}_{\infty},\norm{\D X}_{\infty})$.
Then, using the second assumption in \eqref{3}, we obtain
\begin{align}
\label{Ar1}
P(E_h;B_{\frac 12})-P(\tilde{E}_h;B_{\frac 12})\leq MC\big[\bigl(\varepsilon_h^{\alpha}+o(\varepsilon_h)\bigr)\varepsilon_h
+\varepsilon_h^{1+\alpha}(P(E_h;B_1))^{\alpha}\big].
\end{align}
{We want apply now Theorem \ref{VarP} on the left-hand side. For this reason let us observe that by Lemma \ref{height bound}, for $h$ large enough, $|x_n|<\delta$ for every $x\in \partial E_h$, so that $\psi'=0$ and then we can write 
$$\nabla X(x)=e_n\otimes\nabla'\phi(x'),\quad {\textnormal{div} X=\phi\psi'=0},
$$
thus concluding
$$
\textnormal{div}_{E_h} X=-\langle\nabla X\nu_{E_h},\nu_{E_h}\rangle=
-\langle \nu_{E_h}, e_n\rangle \langle \D'\phi,\nu_{E_h}' 
\rangle \quad
\textnormal{on }
\partial E_h.
$$
Therefore, applying Theorem \ref{VarP}, we obtain
\begin{equation*}
P(E_h;B_{\frac 12})-P(\tilde{E}_h;B_{\frac 12})=
(\varepsilon_h+O(\varepsilon_h^2))\int_{\dd E_h\cap B_{\frac{1}{2}}}\langle \nu_{E_h}, e_n\rangle \langle \D'\phi,\nu_{E_h}' \rangle\,d\mathcal{H}^{n-1},
\end{equation*}
and then inserting this equality in \eqref{Ar1} we deduce,
\begin{align}
\label{Ar2}
(\varepsilon_h+O(\varepsilon_h^2))\int_{\dd E_h\cap B_{\frac{1}{2}}}\langle \nu_{E_h}, e_n\rangle \langle \D'\phi,\nu_{E_h}' \rangle\,d\mathcal{H}^{n-1}\leq MC\big[\bigl(\varepsilon_h^{\alpha}+o(\varepsilon_h)\bigr)\varepsilon_h
+\varepsilon_h^{1+\alpha}(P(E_h;B_1))^{\alpha}\big].
\end{align}
Finally, if we replace $\phi$ by $-\phi$, we deduce dividing by $\varepsilon_h$
$$
\Bigl|\int_{\dd E_h\cap B_{\frac{1}{2}}}\langle \nu_{E_h}, e_n\rangle \langle \D'\phi,\nu_{E_h}' \rangle\,d\mathcal{H}^{n-1}\Bigr|
\leq MC
\bigl(\varepsilon_h^{\alpha}+o(\varepsilon_h)\bigr)\bigl(1+P(E_h;B_1)^\alpha\bigr),
$$
then recalling that $\alpha>\frac{n-1}{n}\geq \frac 12$ we deduce
\begin{equation}
\label{5}
\lim_{h\rightarrow +\infty}\frac{1}{\sqrt{\varepsilon_h}} \Bigl|\int_{\dd E_h\cap B_{\frac{1}{2}}}\langle \nu_{E_h}, e_n\rangle \langle \D'\phi,\nu_{E_h}' \rangle\Bigr|\,d\mathcal{H}^{n-1}=0.
\end{equation}}
Decomposing $\dd E_h\cap B_{\frac{1}{2}}=\big([\Gamma_{f_h}\cup(\dd E_h\setminus\Gamma_{f_h})]\setminus(\Gamma_{f_h}\setminus\dd E_h)\big)\cap B_{\frac{1}{2}}$, we deduce
\begin{equation}
\label{4}
\begin{split}
& -\frac{1}{\sqrt{\varepsilon_h}}\int_{\dd E_h\cap B_{\frac{1}{2}}}\langle \nu_{E_h}, e_n\rangle \langle \D'\phi,\nu_{E_h}' \rangle\,d\mathcal{H}^{n-1}
 =\frac{1}{\sqrt{\varepsilon_h}}\bigg[-\int_{ \Gamma_{f_h}\cap B_{\frac{1}{2}}}\langle \nu_{E_h}, e_n\rangle \langle \D'\phi,\nu_{E_h}' \rangle\,d\mathcal{H}^{n-1}\\
& - \int_{(\dd E_h\setminus \Gamma_{f_h})\cap B_{\frac{1}{2}}}\langle \nu_{E_h}, e_n\rangle \langle \D'\phi,\nu_{E_h}' \rangle\,d\mathcal{H}^{n-1}+\int_{(\Gamma_{f_h}\setminus \dd E_h)\cap B_{\frac{1}{2}}}\langle \nu_{E_h}, e_n\rangle \langle \D'\phi,\nu_{E_h}' \rangle\,d\mathcal{H}^{n-1}\bigg].
\end{split}
\end{equation}
Since by the second inequality in \eqref{1} we have
\begin{equation*}
\bigg |\frac{1}{\sqrt{\varepsilon_h}} \int_{(\dd E_h\setminus \Gamma_{f_h})\cap B_{\frac{1}{2}}}\langle \nu_{E_h}, e_n\rangle \langle \D'\phi,\nu_{E_h}' \rangle\,d\mathcal{H}^{n-1} \bigg | \leq C_8\sqrt{\varepsilon_h}\sup_{\R^{n-1}}\abs{\D'\phi},
\end{equation*}
\begin{equation*}
\bigg |\frac{1}{\sqrt{\varepsilon_h}} \int_{(\Gamma_{f_h}\setminus \dd E_h)\cap B_{\frac{1}{2}}}\langle \nu_{E_h}, e_n\rangle \langle \D'\phi,\nu_{E_h}' \rangle\,d\mathcal{H}^{n-1} \bigg | \leq C_8\sqrt{\varepsilon_h}\sup_{\R^{n-1}}\abs{\D'\phi},
\end{equation*}
then by \eqref{5} and the area formula, we infer
\begin{equation*}
0=\lim_{h\rightarrow +\infty}\frac{-1}{\sqrt{\varepsilon_h}}\int_{ \Gamma_{f_h}\cap B_{\frac{1}{2}}}\langle \nu_{E_h}, e_n\rangle \langle \D'\phi,\nu_{E_h}' \rangle\,d\mathcal{H}^{n-1}=\lim_{h\rightarrow +\infty}\frac{1}{\sqrt{\varepsilon_h}}\int_{ {D_{\frac{1}{2}}}}\frac{\langle \D' f_h,\D'\phi\rangle}{\sqrt{1+\abs{\D' f_h}^2}}\,dx'.
\end{equation*}
This proves that $g$ is harmonic.\\
\textbf{Step 2.} The proof of this step now follows  exactly as in \cite{FJ} using the height bound lemma and the reverse Poincaré inequality. We give here the proof for the sake of completeness.\\ By the mean value property of harmonic functions, Lemma 25.1 in \cite{Ma}, Jensen's inequality, semicontinuity and the third inequality in \eqref{1} we deduce that
\begin{equation*}
\begin{split}
& \lim_{h\rightarrow \infty}\frac{1}{\varepsilon_h}\int_{ {D_{2\tau}}}\abs{f_h(x')-(f_h)_{2\tau}-\langle(\D' f_h)_{2\tau},x'\rangle}^2\,dx'\\
& =\int_{ {D_{2\tau}}}\abs{g(x')-(g)_{2\tau}-\langle(\D' g)_{2\tau},x'\rangle}^2\,dx'\\
& = \int_{ {D_{2\tau}}}\abs{g(x')-g(0)-\langle\D' g(0),x'\rangle}^2\,dx'\\
&\leq c(n)\tau^{n-1}\sup_{x'\in  {D_{2\tau}}}\abs{g(x')-g(0)-\langle\D' g(0),x'\rangle}^2\\
& \leq c(n)\tau^{n+3}\int_{ {D_{\frac{1}{2}}}}\abs{\D' g}^2\,dx'\leq c(n)\tau^{n+3}\liminf_{h\rightarrow \infty}\int_{ {D_{\frac{1}{2}}}}\abs{\D' g_h}^2\,dx'\\
& \leq \tilde{C}(n,C_8)\tau^{n+3}.
\end{split}
\end{equation*}
On one hand, using the area formula, the mean value property, the previous inequality and setting
\begin{equation*}
c_h:=\frac{(f_h)_{2\tau}}{\sqrt{1+\abs{(\D' f_h)_{2\tau}}^2}}, \quad \nu_h:=\frac{(-(\D' f_h)_{2\tau},1)}{\sqrt{1+\abs{(\D' f_h)_{2\tau}}^2}},
\end{equation*}
we have
\begin{equation*}
\begin{split}
&\limsup_{h\rightarrow \infty}\frac{1}{\varepsilon_h}\int_{\dd E_h\cap\Gamma_{f_h}\cap B_{2\tau}}\abs{\langle\nu_h,x\rangle-c_h}^2\,d\mathcal{H}^{n-1}\\
& = \limsup_{h\rightarrow \infty}\frac{1}{\varepsilon_h}\int_{\dd E_h\cap\Gamma_{f_h}\cap B_{2\tau}}\frac{\abs{\langle -(\D' f_h)_{2\tau},x'\rangle+f_h(x')-(f_h)_{2\tau}}}{1+\abs{(\D' f_h)_{2\tau}}^2}^2\sqrt{1+\abs{\D 'f_h(x')}^2}\,dx'\\
& \leq \lim_{h\rightarrow \infty}\frac{1}{\varepsilon_h}\int_{ {D_{2\tau}}}\abs{f_h(x')-(f_h)_{2\tau}-\langle(\D' f_h)_{2\tau},x'\rangle}^2\,dx'\leq \tilde{C}(n,C_8)\tau^{n+3}.
\end{split}
\end{equation*}
On the other hand, arguing as in Step 1, we immediately get from the height bound lemma and the first two inequalities in \eqref{1} that
\begin{equation*}
\lim_{h\rightarrow \infty}\frac{1}{\varepsilon_h}\int_{(\dd E_h\setminus\Gamma_{f_h})\cap B_{2\tau}}\abs{\langle\nu_h,x\rangle-c_h}^2\,d\mathcal{H}^{n-1}=0.
\end{equation*}
Hence we conclude that
\begin{equation}
\label{6}
\limsup_{h\rightarrow \infty}\frac{1}{\varepsilon_h}\int_{\dd E_h\cap B_{2\tau}}\abs{\langle\nu_h,x\rangle-c_h}^2\,d\mathcal{H}^{n-1}\leq \tilde{C}(n,C_8)\tau^{n+3}.
\end{equation}
We claim that the sequence $\lbrace {\mathbf e}_h(0,2\tau,\nu_h) \rbrace_{h\in\N}$ is infinitesimal; indeed, by the definition of excess, Jensen's inequality and the third inequality in \eqref{1} we have
\begin{equation*}
\begin{split}
& \limsup_{h\rightarrow \infty}\int_{\dd E_h\cap B_{2\tau}}\abs{\nu_{E_h}-\nu_h}^2\,d\mathcal{H}^{n-1}\\
& \leq\limsup_{h\rightarrow \infty}\bigg[ 2\int_{\dd E_h\cap B_{2\tau}}\abs{\nu_{E_h}-e_n}^2\,d\mathcal{H}^{n-1}+2\abs{e_n-\nu_h}^2\mathcal{H}^{n-1}(\dd E_h\cap B_{2\tau})\bigg]\\
& \leq \limsup_{h\rightarrow \infty}\bigg[4\varepsilon_h+2\mathcal{H}^{n-1}(\dd E_h\cap B_{2\tau})\frac{\abs{((\D' f_h)_{2\tau},\sqrt{1+\abs{(\D' f_h)_{2\tau}}^2}-1)}^2}{1+\abs{(\D' f_h)_{2\tau}}^2}\bigg]\\
& \leq \limsup_{h\rightarrow \infty}\big[4\varepsilon_h+4\mathcal{H}^{n-1}(\dd E_h\cap B_{2\tau})\abs{(\D' f_h)_{2\tau}}^2\big]\\
& \leq \limsup_{h\rightarrow \infty}\bigg[ 4\varepsilon_h+4\int_{ {D_{\frac{1}{2}}}}\abs{\D' f_h}^2\,dx' \bigg]\leq \lim_{h\rightarrow \infty}[4\varepsilon_h+4C_8\varepsilon_h]=0.
\end{split}
\end{equation*}
Therefore, applying the reverse Poincaré inequality and \eqref{6}, we have for $h$ large that
\begin{equation*}
\begin{split}
{\mathbf e}_h(0,\tau)\leq {\mathbf e}_h(0,\tau,\nu_h)\leq C_9(\tilde{C}\tau^2{\mathbf e}_h(0,1)+\mathcal{D}(0,2\tau)+(2\tau r_h)^\gamma),
\end{split}
\end{equation*}
which is a contradiction if we choose $C_{10}>C_9\max\lbrace\tilde{C},2^\gamma\rbrace$.
\end{proof}

\section{Proof of the main theorem}
\label{Proof of the main theorem}
The proof works exactly as in \cite{FJ}. We give here some details to emphasize the dependence of the constant $\varepsilon$ appearing in the statement of Theorem \ref{Teorema principale} from the structural data of the functional. The proof is divided in four steps.\\
\textbf{Step 1.} We show that for every $\tau\in(0,1)$ there exists $\varepsilon_5=\varepsilon_5(\tau)>0$ such that if ${\mathbf e}(x,r)\leq\varepsilon_6$, then
\begin{equation*}
\mathcal{D}(x,\tau r)\leq C_4\tau\mathcal{D}(x,r),
\end{equation*}
where $C_4$ is from Lemma \ref{Lemma decadimento 1}. Assume by contradiction that for some $\tau\in(0,1)$ there exist two positive sequences $(\varepsilon_h)_h$ and $(r_h)_h$ and a sequence $(E_h,u_h)$ of $(\Lambda r_h^\gamma,\alpha)$-minimizers of $\mathcal{F}_{r_h}$ in $B_1$ with equibounded energies such that, denoting by ${\mathbf e}_h$ the excess of $E_h$ and by $\mathcal{D}_h$ the rescaled Dirichlet integral of $u_h$, we have that $0\in\dd E_h$,
\begin{equation}
\label{Eqn 7}
{\mathbf e}_h(0,1)=\varepsilon_h\rightarrow 0 \quad\text{and}\quad \mathcal{D}_h(0,\tau)>C_4\tau\mathcal{D}_h(0,1).
\end{equation}
Thanks to the energy upper bound (Theorem \ref{Energy upper bound}) and the compactness lemma (Lemma \ref{Lemma compattezza}), we may assume that $E_h\rightarrow E$ in $L^1(B_1)$ and $0\in\dd E$. Since, by lower semicontinuity, the excess of $E$ at 0 is null, $E$ is a half-space in $B_1$, say $H$. In particular, for $h$ large, it holds
\begin{equation*}
|(E_h\Delta H)\cap B_1|<\varepsilon_0(\tau)|B_1|,
\end{equation*}
where $\varepsilon_0$ is from Lemma \ref{Lemma decadimento 1}, which gives a contradiction with the inequality \eqref{Eqn 7}.\\
\textbf{Step 2.} Let $U\subset\subset\Omega$ be an open set. Prove that for every $\tau\in(0,1)$ there exist two positive constants $\varepsilon_6=\varepsilon_6(\tau,U)$ and $C_{11}$ such that if $x_0\in\dd E$, $B_r(x_0)\subset U$ and ${\mathbf e}(x_0,r)+\mathcal{D}(x_0,r)+r^\gamma<\varepsilon_6$, then
\begin{equation}
\label{Eqn 8}
{\mathbf e}(x_0,\tau r)+\mathcal{D}(x_0,\tau r)+(\tau r)^\gamma\leq C_{11}(\tau{\mathbf e}(x_0,r)+\tau\mathcal{D}(x_0,r)+(\tau r)^\gamma).
\end{equation}
Fix $\tau\in(0,1)$ and assume without loss of generality that $\tau<\frac{1}{2}$. We can distinguish two cases. \\
\textit{Case 1:} $\mathcal{D}(x_0,r)+r^\gamma\leq \tau^{-n}{\mathbf e}(x_0,r)$. If ${\mathbf e}(x_0,r)<\min\{ \varepsilon_4 (\tau,\tau^{-n}),\varepsilon_5(2\tau)\}$ it follows from Theorem \ref{Miglioramento eccesso} and Step 1 that
\begin{equation*}
{\mathbf e}(x_0,\tau r)\leq C_{10}(\tau^2{\mathbf e}(x_0,r)+\mathcal{D}(x_0,2\tau r)+(\tau r)^\gamma)\leq C_{10}(\tau{\mathbf e}(x_0,r)+2C_4\tau\mathcal{D}(x_0,r)+(\tau r)^\gamma).
\end{equation*}
\textit{Case 2:} ${\mathbf e}(x_0,r)\leq\tau^{n}(\mathcal{D}(x_0,r)+r^\gamma)$. By the property of the excess at different scales, we infer
\begin{equation*}
{\mathbf e}(x_0,\tau r)\leq \tau^{1-n}{\mathbf e}(x_0,r)\leq (\tau\mathcal{D}(x_0,r)+(\tau r)^\gamma).
\end{equation*}
We conclude that choosing $\varepsilon_6=\min\{ \varepsilon_4(\tau,\tau^{-n}),\varepsilon_5(2\tau),\varepsilon_5(\tau) \}$, inequality \eqref{Eqn 8} is verified.\\
\textbf{Step 3.} Fix $\sigma\in(0,\frac{\gamma}{2})$ and choose $\tau_0\in(0,1)$ such that $C_{11}\tau_0^\gamma\leq\tau_0^{2\sigma}$. Let $U\subset\subset\Omega$ be an open set. We define
\begin{equation*}
\begin{split}
\Gamma\cap U:=\{ x\in\dd E\cap U\,:\, {\mathbf e}(x,r)+\mathcal{D}(x,r)+r^\gamma<\varepsilon_6(\tau_0,U),
\text{ for some }r>0\text{ such that }B_r(x_0)\subset U \}.
\end{split}
\end{equation*}
Note that $\Gamma\cap U$ is relatively open in $\dd E$. We show that $\Gamma\cap U$ is a $C^{1,\sigma}$-hypersurface. Indeed, inequality \eqref{Eqn 8} implies via standard iteration argument that if $x_0\in\Gamma\cap U$ there exist $r_0>0$ and a neighborhood $V$ of $x_0$ such that for every $x\in\dd E\cap V$ it holds:
\begin{equation*}
{\mathbf e}(x,\tau_0^k r_0)+\mathcal{D}(x,\tau_0^k r_0)+(\tau_0^k r_0)^\gamma\leq \tau_0^{2\sigma k}, \quad\text{for }k\in\N_0.
\end{equation*}
In particular ${\mathbf e}(x,\tau_0^k r_0)\leq \tau_0^{2\sigma k}$ and, arguing as in \cite{FJ}, we obtain that for every $x\in\dd E\cap V$ and $0<s<t<r_0$ it holds
\begin{equation*}
|(\nu_E)_s(x)-(\nu_E)_t(x)|\leq ct^\sigma,
\end{equation*}
for some constant $c=c(n,\tau_0,r_0)$, where
\begin{equation*}
(\nu_E)_t(x)=\fint_{\dd E\cap B_t(x)}\nu_E\,d\mathcal{H}^{n-1}.
\end{equation*}
The previous estimate first implies that $\Gamma\cap U$ is $C^1$. By a standard argument we then deduce again from the same estimate that $\Gamma\cap U$ is a $C^{1,\sigma}$-hypersurface. Finally we define $\Gamma:=\cup_i(\Gamma\cap U_i)$, where $(U_i)_i$ is an increasing sequence of open sets such that $U_i\subset\subset\Omega$ and $\Omega=\cup_i U_i$.\\
\textbf{Step 4.} Finally we are in position to prove that there exists $\epsilon>0$ such that
\begin{equation*}
\mathcal{H}^{n-1-\epsilon}(\dd E\setminus \Gamma)=0.
\end{equation*}
Being the argument rather standard, 
Setting $\Sigma=\Big\{ x\in\dd E\setminus\Gamma\,:\, \lim\limits_{r\rightarrow 0}\mathcal{D}(x,r)=0 \Big\}$, by Lemma \ref{Lemma maggiore sommabilità} we have that $\D u\in L^{2s}_{loc}(\Omega)$ for some $s=s(n,\nu,N,L)>1$ and we have that
\begin{equation*}
\text{dim}_{\mathcal{H}}\Big(\Big\{x\in\Omega\,:\, \limsup_{r\rightarrow 0}\mathcal{D}(x,r)>0 \Big\}\Big)\leq n-s.
\end{equation*}
The conclusion follows as in \cite{FJ} (see also \cite{DFR} and \cite{DF0}) showing that $\Sigma=\emptyset$ if $n\leq 7$ and $\text{dim}_\mathcal{H}(\Sigma)\leq n-8$ if $n\geq 8$.


\end{document}